# Cooperative Trajectory Control for Synchronizing the Movement of Two Connected and Autonomous Vehicles Separated in a Mixed Traffic Flow

**Jiahua Qiu**
Department of Civil and Coastal Engineering,
University of Florida, Gainesville, Florida, 32608
Email: jq22@ufl.edu

**Lili Du***
Associate professor
Department of Civil and Coastal Engineering,
University of Florida, Gainesville, Florida, 32608
Email: lilidu@ufl.edu

**ABSTRACT**

When connected and autonomous vehicles (CAVs) are widely used in the future, we can foresee many essential applications, such as platoon formation and autonomous police patrolling, which need two CAVs, originally separated in a mixed traffic flow involving CAVs and human-drive vehicles (HDVs), to quickly approach each other and then keep a stable car-following mode. The entire process should not jeopardize surrounding traffic safety and efficiency. The existing literature has not studied this CAV synchronization control well, and this study seeks to partially make up this gap. To do that, we develop a Model Predictive Control model embedded with a mixed-integer nonlinear program (MINLP-MPC), which integrates micro- and macro-traffic flow models to capture hybrid traffic flow dynamics. Specifically, the MPC will generate control law at each discrete timestamp to manage the microscopic movements of the two subject CAVs while predicting their neighboring vehicles' movement by well-accepted car-following models and estimating the distant upstream traffic' response by the macroscopic traffic model such as cell transmission model (CTM). The MINLP-MPC is multi-objective seeking to sustain both synchronization and traffic efficiencies. To generate such well-balanced optimal control, we noticed that the synchronization experiences two distinct phases, sequentially completing the catch-up and platooning tasks. Accordingly, we transfer MINLP-MPC to a hybrid MPC system consisting of two sequential MPCs, respectively prioritizing the catch-up and platooning control. Then, we develop a weighting strategy to tune the control priorities adaptively. The recursive feasibility of the MPC is mathematically proved. A numerical study built upon the NGSIM dataset demonstrates the efficiency and effectiveness of our approaches under different congestion levels and CAV penetrations.

## 1. Introduction

Connected and Autonomous Vehicles (CAVs), equipped with wireless communication and automation capabilities bring great opportunities to develop cooperative driving (Milanés and Shladover, 2014; Vahidi and Sciarretta, 2018) in the future. Existing literature has investigated this potential for different traffic scenarios, such as vehicle platooning (Čičić and Johansson, 2019; Dey et al., 2015; Liang et al., 2015; Vahidi and Sciarretta, 2018; Wang et al., 2018; Zhang et al., 2021) on a highway segment, cooperative lane change and merging at mainline and on-ramp (or off-ramp) highway intersection (Bevly et al., 2016; Mu et al., 2021; Rios-Torres and Malikopoulos, 2016), and cooperative urban intersection control for CAV fleet (Ashtiani et al., 2018; Yang et al., 2016; Zohdy and Rakha, 2012). This study shares the same research interests in developing the cooperative driving, but particularly focuses on the scenario in which two CAVs want to approach each other quickly and then keep a stable car-following mode, while they are originally apart and separated by other vehicles (e.g., CAVs or human-drive vehicles (HDVs)) in same or different lanes. We name this type of cooperative driving CAV synchronization control (shortened as synchronization control hereafter).

In the era of CAVs in the future, we can foresee many important applications of this synchronization control. The following only exemplifies a few. First, scholars and OEMs have been aware that CAVs from different manufacturers, such as Tesla and Nissan, will most likely install different cruise control algorithms and they are potentially incompatible (Bhoopalam et al., 2018). Then, when a single CAV drives on the road, it may want to approach the other CAV nearby from the same manufacturer to form a short platoon and gain the benefit of the cruise control. Similarly, engineers have noticed that it is not proper to simultaneously manipulate multiple heavy-duty CAVs' motions to form a long truck platoon due to its significant impacts on the safety and efficiency of neighborhood traffic. Instead, an existing study (Liang et al., 2015) considers a practical scheme to have one heavy-duty CAV truck approach another first to form a short platoon. Besides, there are more futuristic applications when CAVs are combined with other emerging technologies. For example, existing studies (Hallmark et al., 2019; Subedi et al., 2020) show that CAVs can potentially be involved in police patrolling activity, where a CAV patrol detects and sends warnings to a speeding CAV and then catches up without triggering a siren to disrupt the traffic. Moreover, by combining CAV and advanced battery-to-battery charging technologies, both scholars and practitioners (Abdolmaleki et al., 2019; Chakraborty et al., 2020; Kosmanos et al., 2018; Maglaras et al., 2014) are investigating emerging mobile vehicle-to-vehicle charging services, which dispatch one electric vehicle (CAV) to approach and then charge another electric vehicle (EV) on-the-move. This new mobile charging service calls for the efficient synchronization control so that we can have one EV on a road catch up the other EV and then maintain a desired relative speed and spacing between the EV pair to sustain cable assembly or stable wireless charging.

On the other hand, this synchronization control presents unique needs which make existing cooperative driving algorithms developed for other traffic scenarios insufficient. More exactly, we consider that the two CAVs initially separated by other vehicles in different lanes want to cooperatively catch up each other first and then keep a stable and efficient car-following mode. Consequently, the synchronization control involves a cooperative catch-up and then platooning control process, each of which introduces complicated lateral and longitudinal movements, e.g., lane-changing and overtaking, etc. Nonetheless, most existing car-following or cooperative lane-charging control (Rios-Torres and Malikopoulos, 2016; Wang et al., 2018) can only handle part of the movement control not the entire synchronization control. Moreover, the entire process of synchronizing two CAVs' movements will affect their neighbor vehicles and even upstream traffic. To account for this traffic impact, both the micro- and macro- traffic dynamics should be well-integrated into the synchronization control. This hybrid traffic effect resulting from CAV movements is not well investigated by existing platoon formation (Čičić and Johansson, 2019; Ganaoui-Mourlan et al., 2021;

Liang et al., 2015; Tuchner and Haddad, 2017) or cooperative lane change and merging (Horowitz et al., 2004; Milanés and Shladover, 2014; Rios-Torres and Malikopoulos, 2016). To the best of our knowledge, this synchronization control has not been well studied in the literature. Only a few studies scratch it but with simple traffic environment setups. Our literature review later presents the research gap further in detail.

This study seeks to bridge the research gap by developing an efficient and optimal synchronization control. It considers that the two subject CAVs are initially separated (i.e., located in a general position) in a mixed flow including both HDVs and CAVs. The synchronization control will generate optimal control laws to guide the two subject CAVs to approach each other smoothly and efficiently without jeopardizing surrounding traffic efficiency and safety. To achieve this research goal, we contribute the following important modeling and mathematical analysis.

First of all, we developed a Model Predictive Control (MPC) model, which repeatedly generates optimal trajectory control laws to guide the movements of the two subject CAVs for this synchronization purpose at each sampling time step according to the initial and predicted system condition in a time horizon. Overall, the MPC involves the enhanced features along with new complexity as follows. It is a hybrid model that captures the microscopic vehicle movement dynamics by existing well-accepted models for the subject CAVs, their neighbor CAVs and HDVs, and then adopts the cell transmission model (CTM) to describe macroscopic traffic flow dynamics. The interaction between the microscopic and macroscopic traffic flow dynamics is captured by introducing integer variables into the CTM so that we can count the flow variation in a cell resulting from microscopic vehicle movements. The resulting hybrid micro-macro model enables the MPC to anticipate the congestion effect resulting from the synchronization control laws but makes some formulations involve nonlinearity. Moreover, the lateral movements of subject CAVs will cause non-smooth car-following relationship changes during the control prediction horizon. To capture this phenomenon, we introduce binary variables to track the leading-following relationships between the subject CAVs and their neighbor vehicles. Together, these enhanced features lead to a highly complicated MPC embedded mixed-integer nonlinear program (MINLP-MPC) with multiple objectives. The complexity of the MINLP-MPC further brings new challenges regarding the recursive feasibility analysis, which is essential to the applicability of an MPC. To address these challenges, we further contribute the following mathematical analysis.

We mathematically proved the recursive feasibility of the MINLP-MPC in Theorem 1 by taking advantage of the problem's unique structure features. It has been known that proving the recursive feasibility of an MPC involving integer variables is difficult in general since it involves enormous state spaces resulting from the combinations of integer variables. We noticed that the integer variables of the MINLP-MPC developed in this study are highly correlated. This observation helped us shrink the state spaces to limited scenarios and then rigorously proved recursive feasibility. More importantly, the proof specifies the conditions that render the recursive feasibility, such as the speed limits, the bound on the human driver's reaction time, and the range of control parameters in the CACC scheme that are followed by neighbor CAVs. The analysis results show that these conditions are satisfied in normal traffic and indicate the strong applicability of the developed MPC system.

Next, the synchronization control is inherently multiple objectives factoring in both microscopic synchronization and macroscopic traffic efficiency, each of which is further interpreted into performance metrics involving several measurements. Accordingly, the optimizer of the MINLP-MPC is multi-objective, requiring effective weighting schemes to generate well-balanced optimal control. To this end, this study contributed an adaptive weighting scheme to balance the traffic efficiency and synchronization efficiency. More exactly, we noticed that the synchronization process would experience two distinct dynamic states: catch-up state and platooning state, each of which has a different control priority in the elements of the performance matrix. Then, we transferred the multi-objective MPC to a hybrid MPC system consisting of two sequential MPCs to decorrelate the objective items and reduce the weighting complexity. Furthermore, we developed a weighting strategy for each MPC system to adaptively tune the control priorities based on the initial traffic condition at each control time step. The effectiveness of this weighting strategy is confirmed in the numerical study.

Last, we carried out extensive numerical experiments to validate the effectiveness of the synchronization control under the developed MPC system. Our experiments show that the adaptive weighting strategy can generate well-balanced control decisions regarding traffic and synchronization performance and outperforms other weighting strategies that have fixed weights. Besides, our numerical experiments demonstrate that the synchronization control is more efficient under a low-moderate congestion level and achieves better performance at high CAV penetration (> 50%).

We present this research by the structure as follows. Following the introduction, Section 2 gives a comprehensive literature review and points out the research gaps we try to bridge. Section 3 provides the preliminary problem formulations, assumptions and briefly introduces the MPC system. Section 4 mathematically formulates the moving dynamics of subject CAVs, HDVs, neighbor CAVs and traffic dynamics of upstream traffic. Next, we formally define the multi-objective MPC system in Section 5 and introduce the adaptive weighting strategy. In Section 6, we analyze the recursive feasibility of the MINLP-MPC model in Section 5. Section 7 further conducts numerical experiments to validate the applicability and effectiveness of our approaches. The entire study and future work are summarized in Section 8.

## 2. Literature review

In recent decades, spurred by the development of CAV technologies, extensive studies have investigated various CAV cooperative driving algorithms for different traffic scenarios. This section will first recognize those relevant studies and then differentiate our study from these existing efforts, including the cooperative adaptive cruise control (CACC) (Badnava et al., 2021), cooperative lane change and merging control (Vahidi and Sciarretta, 2018), and platoon formation (Ganaoui-Mourlan et al., 2021; Liang et al., 2015; Tuchner and Haddad, 2017). Furthermore, considering the MPC integrates micro- and macro-traffic dynamics to account for the comprehensive traffic impact of motion control, we will also review the recent efforts in hybrid model development (Pasquale et al., 2018; Piacentini et al., 2019a, 2019b, 2018). Along with this review, we will highlight the research gaps and challenges which motivate this study.

First of all, the synchronization control is mainly relevant to adaptive cruise control (ACC) and cooperative adaptive cruise control (CACC). ACC controller manages the brake and throttle to control the acceleration of the vehicle based on the driving information of the preceding vehicle. Its capability is further enhanced with the addition of vehicle-to-vehicle (V2V) communication, leading to CACC. CACC exchanges more extensive information with multiple preceding vehicles and cooperatively optimize the driving profile which makes the driving process more stable. Existing studies have developed extensive CACC schemes with various control methods, such as proportional integral derivative control (PID) (Knights et al., 2018), sliding mode control (SMC) (Xiao and Gao, 2011; Yan et al., 2017), and model predictive control (MPC) (Kianfar et al., 2015; Maxim et al., 2017; Negenborn and Maestre, 2014), to regulate spacing and speeds of two adjacent CAVs, aiming to mitigate congestion and energy consumption while keeping the vehicle-string stable (Dey et al., 2015; Vahidi and Sciarretta, 2018; Wang et al., 2018). These controllers are either based on the constant distance (CD) policy (Peters et al., 2014), where the inter-vehicle distance is constant and independent of the vehicle's velocity, or the constant time headway (CTH) policy (Alsuhaim et al., 2021; Naus et al., 2010), where the inter-vehicle distance varies with the vehicle's velocity. It should be noted that the core of an ACC or a CACC system is a vehicle-following control model, which assumes that the two CAVs are well-positioned in a leading and following relationship without other vehicles in-between. Accordingly, they mainly focus on longitudinal movement control and ignore the lateral movement, which is necessary for the synchronization control in this study.

While considering CAVs' lateral movement control, we recognized that many algorithms had been developed to assist cooperative lane change and merging maneuvers. More exactly, optimal lane change decisions can be made in a distributed manner (Nie et al., 2016) or in a centralized framework (Cao et al., 2015), through game models (Talebpour et al., 2015; Wang et al., 2015; Zimmermann et al., 2018), model predictive control approaches(Cao et al., 2015; Liu et al., 2018; Liu and Özgüner, 2015; Zhang et al., 2021), or other methods (Awal et al., 2015; Balal et al., 2016; Choi and Yeo, 2017; Gong and Du, 2016). These studies are mainly developed for two field applications: (1) platoon merging and splitting maneuvers, in which a CAV conducts lane change to join or leave the platoon (Horowitz et al., 2004; Milanés et al., 2010); (2) merging at highway on-ramps, in which a CAV or platoon on-ramp cooperatively merges into a platoon on the mainline road (Rios-Torres and Malikopoulos, 2016). A comprehensive survey can be found in (Rios-Torres and Malikopoulos, 2016; Wang et al., 2018). We noticed that these studies usually either assume a pure CAV or a CAV-and-HDV mixed environment, in which all CAVs, but HDVs are under the control of the lane change or merging control algorithms. Unlike all these related studies, the synchronization control adds one more layer of complexity. It needs to consider at least three types of vehicles mixed in the traffic flow, including the two CAVs following the synchronization control, other CAVs following certain car-following control, and HDVs. Consequently, no existing cooperative lane change or merging control algorithms can address the synchronization control. This new research challenge calls for new methodology development.

The synchronization control is also relevant to the recent studies of platoon formation. Specifically, Tuchner and Haddad (2017) studied platoon formation by optimally controlling the longitudinal and lateral movements of CAVs to form a platoon using interpolating control. However, the study only considers a laboratory environment, and the complexity of the traffic environment is missed. Ganaoui-Mourlan et al. (2021) considered a pure CAV environment and developed a rapidly exploring random trees-based algorithm to generate paths, which organizes CAVs into a platoon. We noticed that the impact of the platoon formation on/from surrounding traffic is ignored. Moreover, several studies (Čičić and Johansson, 2019; Liang et al., 2015) developed catch-up algorithms to have one CAV efficiently approach a leading CAV to form a platoon. However, they assumed the two CAVs were apart away in the same lane without other vehicles in-between. Accordingly, lateral movement control is not considered. In summary, none of the above platoon formation algorithms can address the synchronization control proposed in this study. They only partially accomplish the synchronization control at some steps or under simplified traffic environments.

Last, most existing cooperative driving algorithms/models mainly emphasized CAVs' driving safety and efficiency but ignored their impacts on the macroscopic traffic in a broad range. Having recognized this issue, Zhang et al. (2021) developed a platoon-based cooperative lane-change controller that factors the impact on a long-stretch platoon in the target lane. However, if not properly designed, the synchronization control may inevitably affect CAVs' neighborhoods and upstream traffic safety and efficiency significantly. For example, if the synchronization control is too aggressive, it will make CAVs conduct frequent and improper lane change maneuvers. Consequently, it will lead to significant traffic fluctuations in CAVs' local traffic environment. This negative effect will further generate backward shockwaves and propagate backward to affect upstream traffic. Therefore, the synchronization control must co-consider microscopic and macroscopic traffic flow dynamics. Even though such complicated micro-macro traffic flow dynamics is barely considered by existing studies of CAV cooperative driving, they have been integrated into some traffic flow control (Pasquale et al., 2018; Piacentini et al., 2019a, 2019b, 2018). Specifically, these studies considered a CAV platoon as a moving bottleneck (MB). Then, by manipulating

CAVs' traffic speed, they imposed variable speed limit (VSL) schemes on other vehicles to mitigate macroscopic congestion. To capture macroscopic traffic response to the VSL, the studies of Piacentini et al., (2018, 2019a, 2019b) introduced an Ordinary Differential Equation (ODE) to couple the MBs' longitudinal movement with LWR model. To be noted, the proposed synchronization control involves lane-changing movements. Therefore, we need more new approaches to capture the impacts of CAVs' longitudinal and also lateral movements on macroscopic traffic flow.

In short, the synchronization control enables some particular tasks for which we want to efficiently bring two CAVs close to each other and well position them in a stable leading-following mode. There are many such applications when transportation systems widely use CAVs in the future. Those unique features make it different from most existing traffic scenarios that have been investigated by the current CAV cooperative driving control. This study is motivated by viewing the above application needs and research gaps. Below we provide our methodology development and analysis in detail.

## 3. Problem Formulation

Formally, this study considers two subject CAVs are driving in a mixed traffic flow involving both HDVs and CAVs on a freeway segment between two intersections or two ramps. Figure 1 shows an example. While HDVs follow adjacent vehicles according to humans' response to nearby traffic, CAVs are under certain cruise control such as the CACC scheme. These two subject CAVs, once within a communication range, will trigger a synchronization control to complete a special task, which requires the following subject CAV to approach the leading subject CAV and then pair to a short platoon efficiently. The special task can be on-the-move EV energy transfer (Abdolmaleki et al., 2019), patrol activities (Hallmark et al., 2019; Subedi et al., 2020), or early period of platoon formation as we mentioned in the introduction. The synchronization control seeks to generate optimal control laws for the two subject CAVs so that they can synchronize their movements smoothly and efficiently without jeopardizing surrounding traffic's efficiency and safety. One of the key challenges to the development of such synchronization control is modeling the interaction between the subject CAVs and the surrounding traffic. More exactly, we cannot control every vehicle for this special task. The surrounding traffic is dynamic and involves various uncertainty. To make the problem tractable, we scope the surrounding traffic of the two subject CAVs from two aspects: (i) a vehicle fleet consisting of neighboring HDVs and CAVs within the subject CAVs' monitoring range, for which we can predict their microscopic movements, and (ii) an adjacent upstream macroscopic traffic segment, for which we can only estimate the traffic propagation. Built upon that, we develop a model predictive control (MPC), which generates optimal discrete trajectory control for the two subject CAVs to efficiently pair to a platoon without worsening the hybrid traffic environment. To mathematically formulate this problem, we introduce the following notations and specifications.

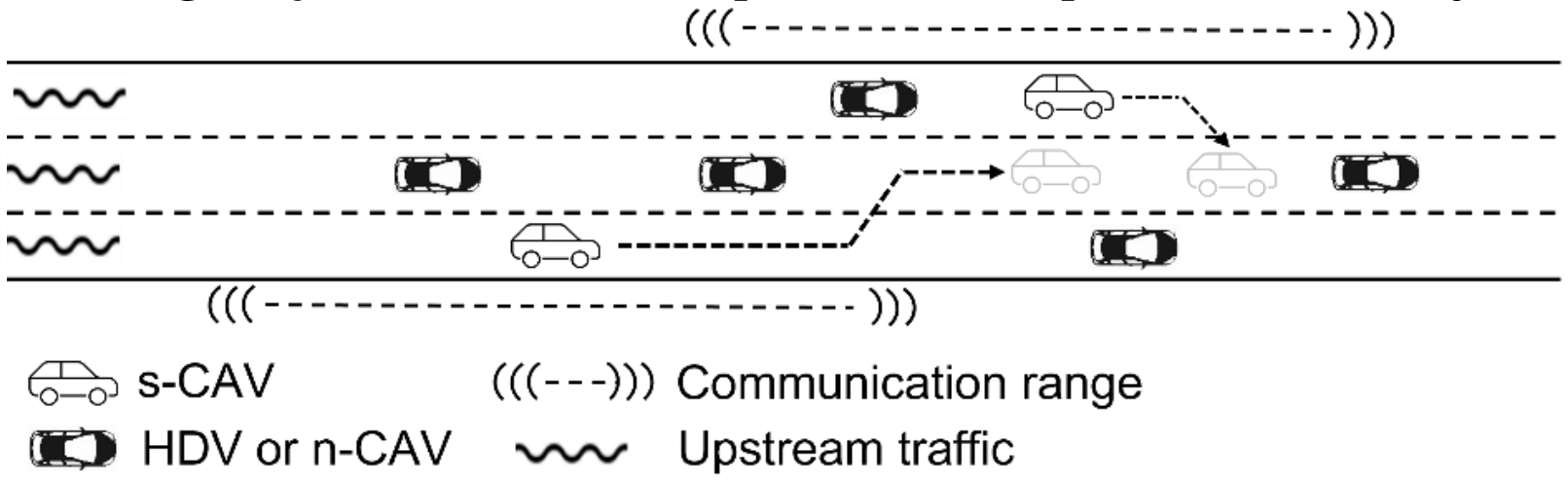

**Figure 1.** pairing two CAVs to a platoon in a mixed traffic environment.

We specify a straight roadway with $L$ number of lanes. The MPC will generate the trajectory control law for the two subject CAV at discrete sampling time steps $t \in \mathbb{Z}_+ \coloneqq \{0,1,2,\dots\}$ with a uniform time interval as $\delta > 0$. We label the two subject CAVs by the indices, $i \in I = \{1,2\}$, where $i = 1$ or 2 represents a leading or following subject CAV. Vectors $z_i(t), v_i(t), u_i(t)$ denote the position, speed, and control input (acceleration/deceleration) of the subject CAV $i$ in longitudinal $x$ and latitudinal $y$ directions at time step $t$, respectively. Namely, $z_i(t) = [x_i(t), y_i(t)]^T, v_i(t) = \left[v_i^x(t), v_i^y(t)\right]^T$, and $u_i = [u_i^x(t), u_i^y(t)]^T, i = 1,2$. The two subject CAVs are able to monitor the movements of their neighboring vehicles $\omega \in \Omega$ within the communication range, which can be either HDVs defined as the set $\Omega_h$ or CAVs defined as the set $\Omega_c$ running in different lanes (we have $\Omega = \{\Omega_c, \Omega_h\} = \left\{\Omega_c^l, \Omega_h^l\right\}_{l=1}^{L}$). Specifically, we denote $\hat{x}_\omega(t)$ and $\hat{v}_\omega(t)$ as the position and speed of a neighboring vehicle $\omega$, $\forall \omega \in \Omega$. To differentiate subject CAVs and their neighboring CAVs, we denote a subject CAV as $s$-CAV, and a CAV in the neighborhood of the subject CAVs as $n$-CAV throughout the rest of the paper. We also consider a vehicle set $J = I \cup \Omega$, which consists of all $s$-CAVs and their neighbor vehicles ($n$-CAVs and neighbor HDVs).

Given that the latitudinal movements of $s$-CAVs will change their lane positions and then affect the following behavior of their neighboring vehicles, we thus introduce binary variables $\gamma_j^l(t), \rho_{j,j'}^l(t)$ and $\eta_{j,j'}^l(t)$ in (1) - (3) to capture the state change resulting from this vehicle movement dynamics.

$$\gamma_j^l(t) = \begin{cases} 1 & \text{if vehicle j run in lane l at time step t} \\ 0 & \text{otherwise} \end{cases}, \qquad \forall j \in J, l \in L \tag{1}$$

$$\rho_{j,j'}^l(t) = \begin{cases} 1 & \text{if vehicle j' is downstream to vehicle j in lane l at time step t} \\ 0 & \text{otherwise} \end{cases}, \forall j, j' \in J, l \in L \tag{2}$$

$$\eta_{j,j'}^l(t) = \begin{cases} 1 & \text{if vehicle j follows j' in lane l at time step t} \\ 0 & \text{otherwise} \end{cases}, \qquad \forall j, j' \in J, l \in L \tag{3}$$

Mainly, $\gamma_j^l(t)$ indicates if a vehicle $j$ is running in a lane $l$ at a time step $t$. $\rho_{j,j'}^l(t)$ demonstrates if vehicle $j'$ is running downstream to vehicle $j$ and $\eta_{j,j'}^l(t)$ specifies if vehicle $j$ is an immediate follower of vehicle $j'$. We determine the variable $\gamma_j^l(t)$ based on the latitudinal coordinate $y_j(t)$ of vehicle $j$ by the logical condition (4), where $[\underline{y}^l, \bar{y}^l]$ is the latitudinal boundary coordinates of lane $l$.

$$\gamma_j^l(t) = [\underline{y}^l < y_j(t) \leq \bar{y}^l], \qquad \forall j, l \in L \tag{4}$$

Clearly, there exists a logical relationship between $\gamma_j^l(t)$ and $\rho_{j,j'}^l(t), \forall j, j' \in J$. For example, $\rho_{j,j'}^l(t) = 1$ if and only if both $j$ and $j'$ are in the same lane, i.e., $\gamma_j^l(t) = \gamma_{j'}^l(t) = 1$, and vehicle $j'$ is in the downstream of vehicle $j$, i.e., $x_j(t) < x_{j'}(t)$. This logical condition is represented in (5), where $\wedge$ represents the logical relationship: "and". The logical conditions (5) can be molded as a set of integer constraints in the Appendix A.

$$\rho_{j,j'}^l(t) = \left[x_j(t) < x_{j'}(t)\right] \wedge \gamma_j^l(t) \wedge \gamma_{j'}^l(t), \qquad \forall j, j' \in J, l \in L \tag{5}$$

Built upon $\rho_{j,j'}^l(t)$, we can determine the value of $\eta_{j,j'}^l(t)$ based on the logical condition in (6). Mainly, we can tell that vehicle $j$ is following vehicle $j'$ at time $t$ on lane $l$ (i.e., $\eta_{j,j'}^l(t) = 1$), if and only if two conditions are satisfied: (i) vehicle $j'$ is in the downstream of vehicle $j$ ($\rho_{j,j'}^l(t) = 1$) and (ii) there are no vehicles between them (i.e., $\left[\sim\rho_{j,k}^l(t)\right] \wedge \left[\sim\rho_{k,j'}^l(t)\right], \forall k \neq j$ and $k \neq j'$). Here, the symbol "~" represents the "not" operation. This logic relationship can be further modeled by integer constraints given in the Appendix A.

$$\eta_{j,j'}^{l}(t) = \rho_{j,j'}^{l}(t) \wedge \cap_{k \in J/\{j,j'\}} \left\{ \left[ \sim \rho_{j,k}^{l}(t) \right] \wedge \left[ \sim \rho_{k,j'}^{l}(t) \right] \right\}, \qquad \forall j,j' \in J, l \in L \quad (6)$$

Besides the $s$-CAVs and their neighboring vehicles (HDVs and $n$-CAVs), the synchronization control also considers the flow dynamics of the upstream traffic outside of the monitoring range of $s$-CAVs. To do that, we discretize the highway segment without intermediate entrances or exits into $C$ number of cells for each lane $l \in L$. Each cell has length $\Delta L$. According to Courant-Friedrichs-Lewy condition, the conservation law holds only when $\Delta L \geq \delta v_f$ (Daganzo, 1995). Each cell $c = 1,2,\ldots,C$ at a time step $t = 1,2,\ldots T$ is associated with the traffic characteristics as follows. $\mathscr{f}_c^l(t)$: the traffic flow of cell c in lane l [veh/h]; $\mathscr{y}_c^l(t)$: entering flow to cell $c$ from cell $c-1$ in lane $l$ [veh/h]; $\mathscr{s}_c^l(t)$: sending flow from $c$ to cell $c+1$ [veh/h] in lane $l$; $\mathscr{r}_c^l(t)$: receiving flow of cell $c$ from cell $c-1$ in lane $l$ (receiving flow) [veh/h].

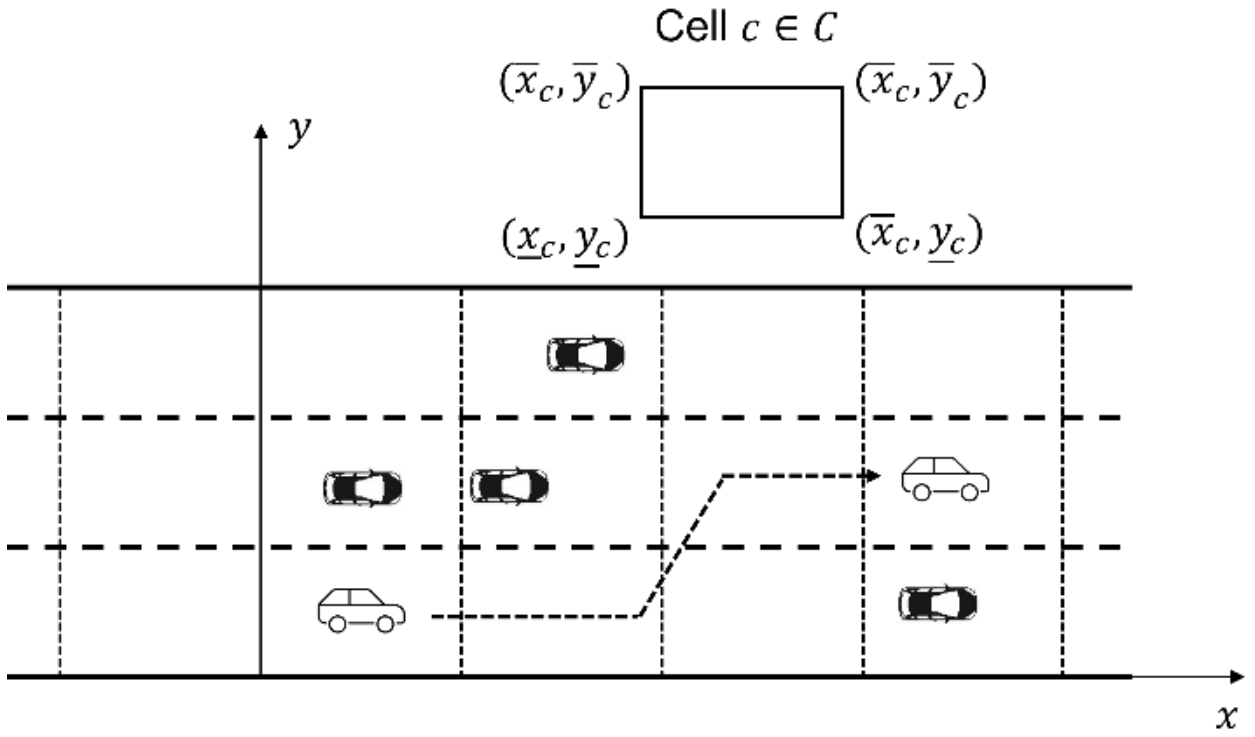

Figure 2 roadway representation in CTM.

The sections below will introduce the microscopic vehicle dynamics, the macroscopic traffic dynamics, and their interaction. Built upon that, we developed the MPC for this synchronization control.

## 4. Microscopic and Macroscopic Traffic Dynamics

The synchronization control will generate the trajectory control law for the $s$-CAVs while predicting the moving dynamics of neighbor vehicles, i.e., $n$-CAVs and neighboring HDVs, as well as the upstream traffic propagation. This section develops mathematical formulations to capture the corresponding microscopic and macroscopic traffic dynamics. To be noted, the latitudinal movements of $s$-CAVs enable lane-change behaviors. It raises the difficulties in formulating the traffic dynamics of its neighboring vehicles and upstream traffic. We address those challenges in the following subsections.

### 4.1 Subject CAVs dynamics

The movements of the $s$-CAVs are subject to the longitudinal and latitudinal dynamics as well as several important states and control constraints. More exactly, we describe the longitudinal and latitudinal dynamics of the $s$-CAVs by the well-accepted double-integrator model in (7) and (8) (Nilsson et al., 2016, 2015; Tomas-Gabarron et al., 2013; Wang et al., 2019).

$$z_i(t+1) = z_i(t) + \delta v_i(t) + \frac{\delta^2}{2} u_i(t), \qquad \forall i \in I \quad (7)$$

$$v_i(t+1) = v_i(t) + \delta u_i(t), \qquad \forall i \in I \quad (8)$$

where $z_i(t), v_i(t), u_i(t)$ represents $z_i(t\delta), v_i(t\delta), u_i(t\delta)$ respectively for notational simplicity. Moreover, we consider that the control input $u_i(t)$ and speed $v_i(t)$ keep constant in each time interval $[t, t+1)$ for $t \in T$, and they are subject to the control and speed limits in (9) and (10).

$$\underline{a} \le u_i(t) \le \overline{a}, \qquad \forall i \in I \tag{9}$$

$$\underline{v} \le v_i(t) \le \overline{v}, \qquad \forall i \in I \tag{10}$$

Besides, a subject CAV's movement is also subject to the safety distance constraint in (11).

$$x_j(t) - x_i(t) \ge L_i + \tau_i v_i^x(t) - \frac{\left[v_i^x(t) - \underline{v_x}\right]^2}{2\,\underline{a_x}}, \qquad \forall i \in I \tag{11}$$

where $x_j(t)$ is the longitudinal position of the leading vehicle $j$ of subject CAV $i$ follows, $\tau_i > 0$ is the constant reaction time due to dynamic and other delays, and $L_i > 0$ is a constant depending on vehicle length. Constraint (11) is a conservative safety constraint. It specifies a safety lower bound on the safety spacing and ensures conflict-free driving (Gong et al., 2016; Gong and Du, 2018). This safety is suitable for freeway traffic but might be too conservative for urban traffic. It can be adapted for urban traffic by using other safety constraints such as (Wang et al., 2020, 2014; Wu et al., 2020). It won't affect the applicability of our approach.

Constraint (11) indicates the need to track the longitudinal position (i.e., $x_j(t)$) of s-CAV ($i$)'s leading vehicle $j$ to keep a safe distance. However, the leading vehicle of a $s$-CAV can be another $s$-CAV, a $n$-CAV, or an HDV and it may switch during the control prediction zone $T$ due to its latitudinal movement. To capture these dynamics, we utilize binary variables $\gamma_i^l(t)$ and $\rho_{i,j}^l(t)$ as defined in and (1) and (2). Specifically, variable $\gamma_i^l(t)$ is introduced to track the lane that the $s$-CAV($i$) is running in, and variable $\rho_{i,j}^l(t)$ monitors its downstream neighbor vehicles. Consequently, we rewrite the safety constraint (11) into constraints (12) and (13), considering two scenarios: (1) $s$-CAV($i$) follows another $s$-CAV($i'$) and (2) $s$-CAV($i$) follows a $n$-CAV or an HDV.

$$x_{i'}(t) - x_i(t) \ge L_i + \tau_i v_i^x(t) - \frac{\left[v_i^x(t) - \underline{v_x}\right]^2}{2\,\underline{a_x}} - M[3 - \gamma_i^l(t) - \gamma_{i'}^l(t) - \rho_{i,i'}^l(t)], \quad \forall i, i' \in I, l \in L \tag{12}$$

$$\hat{x}_{\varpi}(t) - x_i(t) \ge L_i + \tau_i v_i^x(t) - \frac{\left[v_i^x(t) - \underline{v_x}\right]^2}{2\,\underline{a_x}} - M[2 - \gamma_i^l(t) - \rho_{i,\varpi}^l(t)], \forall \varpi \in \Omega_h^l, l \in L, i \in I \tag{13}$$

where, $\hat{x}_{\varpi}(t)$ in constraint (13) is the longitudinal coordinate of neighboring vehicle $\varpi$ at time interval $t$. $M$ takes a relatively big value, i.e., the length of the roadway. Constraint (12) is activated when $\gamma_i^l(t) = \gamma_{i'}^l(t) = \rho_{i,i'}^l(t) = 1$, indicating $s$-CAV($i$) follows another $s$-CAV($i'$) in the downstream. Constraint (13) is activated when $\gamma_i^l(t) = \rho_{i,\varpi}^l(t) = 1$, indicating when $s$-CAV($i$) follows a neighbor vehicle $\varpi$ (HDV or $n$-CAV) in the downstream. The two constraints are inactive when $s$-CAV($i'$) and neighboring vehicle $\varpi$ is not in the downstream of $s$-CAV($i$).

### 4.2 HDVs moving dynamics

We consider that neighboring HDVs are not under a trajectory control scheme or our synchronization control. Still, their movements are restrained by their leading vehicles in this local area, thus can be predicted by using Newell's car-following model (14) during the prediction horizon.

$$\hat{v}_{\omega}(t+k) = \frac{\delta}{\tau_{\omega}} s_{\omega,j}(t) - \frac{\delta d_{\omega}}{\tau_{\omega}} \qquad \forall \omega \in \Omega_h^l, l \in L, j \in J \tag{14}$$

where $s_{\omega,j}(t) = x_j(t) - \hat{x}_{\omega}(t), j \in J$ represents the spacing with the HDV $\omega$ and its leading vehicles $j$ at time step $t$; $\tau_{\omega}$ represents the reaction time and $d_{\omega}$ represents the minimum stop distance of HDV $\varpi$. Both $\tau_\omega$ and $d_\omega$ can be adaptively learned by adopting an online-learning approach (Yu et al., 2019). $\tau_\omega$ will be rounded to an integer number, and the control step size $\delta$ will be selected such that $k = \frac{\tau_\omega}{\delta}$ is an integer to fit HDV dynamics into the model.

Please note that (14) does not explicitly identify the HDV's lane. It will be specified when we recognize its leading vehicle $j$ through the logic constraints in (4). Mainly, the leading vehicle $j$ can also be a $s$-CAV ($i$), $n$-CAVs, or HDV $\omega'$. Accordingly, we set $s_{\omega,i}(t) = x_i(t) - \hat{x}_\omega(t)$, if HDV $\omega$ is following a $s$-CAV($i$) in the same lane $l$ (i.e., $\rho^l_{\omega,i}(t) = 1$) at time $t$; otherwise, we make $s_{\omega,i}(t) = M$, which is a sufficiently large number. Similarly, we have $s_{\omega,\omega'}(t) = \hat{x}_{\omega'}(t) - \hat{x}_\omega(t)$, if HDV $\omega$ is following a neighbor vehicle $\omega'$ on the same lane at time step $t$; otherwise, we make $s_{\omega,\omega'}(t) = M$. We formulate the above thoughts in (15).

$$s_{\omega,j}(t) = \begin{cases} x_j(t) - \hat{x}_\omega(t) & \rho^l_{\omega,j}(t) = 1 \\ M & \rho^l_{\omega,j}(t) = 0 \end{cases}, \qquad \forall \omega \in \Omega^l_h, j \in I \cup \Omega^l, l \in L \tag{15}$$

where $\rho^l_{\omega,j}$ is defined by (2) and is determined by the logical condition (5). In particular, (15) considers $\rho^l_{\omega,i}(t)$ regarding $s$-CAV($i$) as a variable because the latitudinal movement of $s$-CAVs will change this following relationship during the control prediction horizon. But, (15) keeps $\rho^l_{\omega,\omega'}(t)$ regarding neighboring HDV($\omega'$) as a given parameter according to the initial condition since we assume that the neighbor HDVs will stay in the same lane during a prediction horizon (see the justification later). Then we identify the following spacing, $s_{\omega,j}(t)$ of HDV $\omega$ by Equation (16), where $s_{\omega,1}(t), s_{\omega,2}(t)$ are spacing between HDV $\omega$ and $s$-CAV(1) and $s$-CAV(2).

$$s_{\omega,j}(t) = min\{s_{\omega,\omega'}(t), s_{\omega,1}(t), s_{\omega,2}(t)\} \qquad \forall \omega \in \Omega^l_h, l \in L \tag{16}$$

Based on (14) to (16), Equations (17) and (19) estimate the trajectory profiles of HDVs based on Newell's car-following model. Please note that the Newell's car-following model does not set speed limits as it is mostly applied in congested scenarios. In this study, when the traffic is sparse and the spacing $s_\omega(t)$ is large, we consider an upper bound $\overline{v_x}$, i.e., free-flow speed, to the estimated speed to make it consistent with reality. Also, we consider a lower bound $\underline{v_x}$ to the estimated speed.

$$\hat{v}_\omega(t+k) = min\left\{\frac{\delta}{\tau_\omega} s_{\omega,j}(t) - \frac{\delta d_\omega}{\tau_\omega}, \overline{v_x}\right\} \qquad \forall \omega \in \Omega^l_h, l \in L \tag{17}$$

$$\hat{v}_\omega(t+k) = max\left\{\frac{\delta}{\tau_\omega} s_{\omega,j}(t) - \frac{\delta d_\omega}{\tau_\omega}, \underline{v_x}\right\} \qquad \forall \omega \in \Omega^l_h, l \in L \tag{18}$$

$$\hat{x}_\omega(t+k) = \hat{x}_\omega(t) + \frac{\tau_\omega}{\delta} \hat{v}_\omega(t), \qquad \forall \omega \in \Omega^l_h, l \in L \tag{19}$$

Finally, HDVs are also subject to the safety distance constraint when the $s$-CAV cuts in front of the HDV. Let us consider a scenario where the $s$-CAV($i$) travels at the speed $v^x_i$ and cuts in front of the HDV who travels faster at the speed $\hat{v}_\omega$, i.e., $\hat{v}_\omega > v^x_i$. Given that the HDV needs time $\tau_\omega$ to react to the lane change of $s$-CAV($i$), it is necessary to specify the lane-change conditions for $s$-CAV($i$) to avoid the collision. Specifically, we require $s$-CAVs to satisfy the following conditions to cut in front of an HDV. They are derived from the safety lane-change constraints in (Chen et al., 2013) (see Section 6 for details).

$$x_i(t) - \hat{x}_\omega(t) \geq d_\omega + \frac{\tau_\omega}{\delta} \hat{v}_\omega(t) - M\left(1 - \eta^l_{\omega,i}(t) + \eta^l_{\omega,i}(t-1)\right) \qquad \forall \omega \in \Omega^l_h, l \in L, i \in I \tag{20}$$

$$v_i(t) - \hat{v}_\omega(t) \geq -M\left(1 - \eta^l_{\omega,i}(t) + \eta^l_{\omega,i}(t-1)\right) \qquad \forall \omega \in \Omega^l_h, l \in L, i \in I \tag{21}$$

where $\eta^l_{\omega,i}(t)$ is defined by (3) and is determined by the logical condition (6). In particular, when the $s$-CAV($i$) cuts in front of HDV($\omega$) at time step $t$ ($\eta^l_{\omega,i}(t) = 1$ and $\eta^l_{\omega,i}(t-1) = 0$), constraint (20) is activated such that the spacing between $s$-CAV($i$) and HDV($\omega$) needs to no smaller than the safety distance $d_\omega + \tau_\omega \hat{v}_\omega(t)$. Also, the speed of $s$-CAV($i$) needs to no less than that of HDV($\omega$) indicated by constraint (21).

Note that the assumption that the neighbor HDVs stay in the same lane during the prediction horizon is valid for our MPC introduced in the later section. Briefly, MPC is a feedback control. At each control sampling time step, it takes the new initial condition of the system as input and generates the control law for the system in multiple steps but only implements the first step. Therefore, if a lane change occurs during a prediction horizon, the correction will be built in the next time step as we recompute the control law by resetting the initial conditions in the MPC. This assumption also makes Equation (16) valid because $\omega'$ is defined as the leading vehicle of HDV $\omega$ at beginning, during the prediction horizon, the leading vehicle of $\omega$ can only be $\omega'$ or one of the subject CAVs.[1]

### 4.3 Neighboring CAVs moving dynamics

Recall that we consider the neighboring CAVs follow certain cruise control schemes, such as CACC. There are many CACC schemes developed in literature (Rahman et al., 2017; Sharath and Velaga, 2020) and we apply a first-order model derived by PATH to model the CACC car-following behavior (Milanés and Shladover, 2014; Qin et al., 2021), since existing field experiments have confirmed its good performance to capture the car-following movement of current commercially available CACC systems. Specifically, the CACC model formulates the vehicle speed by the speed and gap error in the previous time step as follows,

$$\hat{v}_\omega(t+1) = \hat{v}_\omega(t) + k_1 e_\omega(t+1) + k_2 \dot{e}_\omega(t+1) \qquad \forall \omega \in \Omega_c^l, l \in L \qquad (22)$$

where $\omega$ is the neighbor CAV and control gains $k_1$ and $k_2$ are set as 0.01 $s^{-2}$ and 1.6 $s^{-1}$ (Xiao et al., 2017); $\dot{e}_\omega(t)$ is the derivative of the gap error defined in (24).

$$e_\omega(t) = \hat{x}_{\omega'}(t) - \hat{x}_\omega(t) - t_d \hat{v}_\omega(t) \qquad \forall \omega \in \Omega_c^l, l \in L \qquad (23)$$

The derivative of gap error is defined in (24).

$$\dot{e}_\omega(t) = \hat{v}_j(t) - \hat{v}_\omega(t) - t_d a_\omega(t) \qquad \forall \omega \in \Omega_c^l, l \in L \qquad (24)$$

where vehicle $j$ is the leading vehicle of the neighbor CAV $\omega$, $t_d$ is the desired time gap of the CACC controller, which is set to 0.6$s$; $a_\omega(t)$ is the acceleration of CAV $\omega$ at time step $t$, which is represented by $\frac{\hat{v}_\omega(t+1)-\hat{v}_\omega(t)}{\Delta t}$. Combined Equation (22), (24) and (24), we represent the speed of CAV $\omega$ by the speed of $\varpi$ and leading vehicle $j$ at previous time step and its spacing with the leading vehicle, $s_{\omega,j}(t)$.

$$\hat{v}_\omega(t+1) = A\hat{v}_\omega(t) + B\hat{v}_j(t+1) + Cs_{\omega,j}(t+1) \qquad \forall \omega \in \Omega_c^l, l \in L, j \in J \qquad (25)$$

where $A = \frac{\delta(1-k_1 t_d - k_2) + k_2 t_d}{\delta + k_2 t_d}$, $B = \frac{\delta k_2}{\delta + k_2 t_d}$, $C = \frac{\delta k_1}{\delta + k_2 t_d}$.

Again, due to the lateral movement of $s$-CAVs, the leading vehicle of a $n$-CAV can also be a $s$-CAV, an HDV or another $n$-CAV, and it may vary at each time step during the control horizon. Then to implement Equation (25) in our MPC, we also need to recognize the leading vehicle and their trajectory for each $n$-CAV at each time step. Thus, we utilize binary variables, $\eta_{\omega,\omega'}^l(t)$, $\eta_{\omega,1}^l(t)$ and $\eta_{\omega,2}^l(t)$ defined in (3) to help. Specifically, $\eta_{\omega,\omega'}^l(t) = 1$, if the leading vehicle of $n$-CAV($\omega$) is vehicle $\omega'$ (HDV or $n$-CAV);

[1] If the leading vehicle of $\omega$ is a $s$-CAV at beginning of prediction horizon, then $\omega'$ is defined as the nearest HDV or the $n$-CAV in the downstream.

otherwise $\eta^l_{\omega,\omega'}(t) = 0$. $\eta^l_{\omega,i}(t) = 1$, if the leading vehicle of $n$-CAV($\omega$) is a $s$-CAV at time step $t$; otherwise $\eta^l_{\omega,i}(t) = 0$. Consequently, Equation (25) is reformulated to (26) for adapting it to this study.

$$\hat{v}_\omega(t+1) = A\hat{v}_\omega(t) + B\eta^l_{\omega,\omega'}(t)\hat{v}_{\omega'}(t) + B\sum_{i\in I}\eta^l_{\omega,i}\, v^x_i(t) + Cs_\omega(t), \forall\omega \in \Omega^l_c, l \in L \quad (26)$$

where, $s_{\omega,j}(t)$ is captured by Equation (16) also, but only for $\omega \in \Omega^l_c$ and $\eta^l_{\omega,j}(t)$ is identified by the logical condition defined in (6).

Finally, we also require $s$-CAV to satisfy the following lane-change condition (27) to cut in front of an $n$-CAV without collision (see Section 6 for a detailed derivation). Note that most of the safety lane-change constraints are designed for the lane-change vehicle that cuts in front of an HDV, which can not be applied to this scenario where the following vehicle in the target lane is another CAV following a CACC control.

$$v^x_i(t) - \hat{v}_\omega(t) \geq \frac{\tilde{a}_i(t-1)\delta(B-\frac{1}{2})}{A} - M\left(1 - \eta^l_{\omega,i}(t) + \eta^l_{\omega,i}(t-1)\right) \quad \forall\omega \in \Omega^l_c, l \in L, i \in I \quad (27)$$

where $\tilde{a}_i(t-1) = max\left[\underline{a_x}, \underline{a_{i,v}}(t-1)\right]$, and $\underline{a_{i,v}}(t-1) = \frac{\underline{v_x} - v^x_i(t-1)}{\delta}$. When $s$-CAV cuts in front of the $n$-CAV, i.e., $\eta^l_{\omega,i}(t) = 1$ and $\eta^l_{\omega,i}(t-1) = 0$, constraint (27) is activated and ensures that the $s$-CAV should cut-in with enough speed to avoid the collision.

**4.4 Flow dynamics**

The synchronization control will instruct subject CAVs to speed up or slow down, and perform overtaking or lane-change maneuvers so that these two CAVs can approach each other efficiently. However, these maneuvers may harm upstream surrounding traffic. To capture this effect, we will extend standard CTM to describe traffic propagation corresponding to the synchronization control, namely, the interaction between the macroscopic upstream traffic and microscopic subject CAVs/neighbor vehicles' movement. For completeness, we present the CTM in constraints (28) - (31), which capture flow propagation among cells in general. Note that we develop CTM on each lane on the road segment of interest.

$$\mathcal{f}^l_c(t+1) = \mathcal{f}^l_c(t) + \mathcal{y}^l_c(t) - \mathcal{y}^l_{c+1}(t), \quad \forall c \in C_l, l \in L \quad (28)$$

$$\mathcal{y}^l_c(t) = \min\{\mathcal{s}^l_{c-1}(t), \mathcal{r}^l_c(t)\}, \quad \forall c \in C_l, l \in L \quad (29)$$

$$\mathcal{s}^l_c(t) = \min\{\mathcal{f}^l_c(t), q^{l,max}_c\}, \quad \forall c \in C_l, l \in L \quad (30)$$

$$\mathcal{r}^l_c(t) = \min\left\{w^l_c\left(k^{l,max}_c - k^l_c(t)\right), q^{l,max}_c\right\}, \quad \forall c \in C_l, l \in L \quad (31)$$

where $q^{l,max}_c$ is the capacity of the cell $c$ in lane $l$, $w^l_c$ is the congestion wave speed and $k^{l,max}_c$ is the jam density. $k^l_c(t)$ is the density of cell $c$ in lane $l$ on time interval $t$. Note that we develop CTM on each individual lane on the road segment of interest.

For a cell $c$ at time step $t$, constraints (28) is the flow conservation constraint. Namely, the cell $c$'s flow at time step $t+1$, $\mathcal{f}^l_c(t+1)$, is equal to the cell's flow at time step $t$, $\mathcal{f}^l_c(t)$, plus the flow propagation occurring during the time interval from $t$ to $t+1$ (i.e., plus the flow entering to this cell from an upstream cell, $\mathcal{y}^l_c(t)$, and minus the flow leaving from this cell for its downstream cell, $\mathcal{y}^l_{c+1}(t)$ at time step $t$). Constraint (29) limits that at time step $t$, the entering flow to cell $c$ from an upstream cell $c-1$ is given as the minimum value of the sending flow from its upstream cell $c-1$, $\mathcal{s}^l_{c-1}(t)$, and its receiving flow $\mathcal{r}^l_c(t)$ at time step $t$. For a given cell $c$ at a time step $t$, the sending flow $\mathcal{s}^l_c(t)$ is limited by the minimum flow that the cell has and its flow capacity. Constraint (31) further models the receiving flow $\mathcal{r}^l_c(t)$, which is the minimum flow of the cell $c$' capacity and its rest flow capacity at time step $t$.

Constraints (28) - (31) altogether describe the upstream macroscopic traffic flow dynamics. However, this study is particularly interested in capturing the influence of the individual vehicles under the synchronization control or monitoring, i.e., $s$-CAVs ($i \in I$) and their neighbor vehicles ($\omega \in \Omega$), on the upstream flow propagation. Specifically, the interactions between those individual vehicles and macroscopic upstream traffic first occur in boundary cells where they are mixed. For these cells, the entering (or leaving) of those individual vehicles ( $i \in I$ or $\omega \in \Omega$) increases (or reduces) the rest flow capacity of the cell and limits the upstream flow propagation from upstream cells. For example, when the deceleration of a subject CAV slows down its following neighbor vehicles in a cell $c$, the corresponding density, $k_c^l(t)$ will increases. The increased cell density will reduce the receiving flow $r_c^l(t)$ from its upstream cell according to constraint (31), which further reduces the entering flow to the cell $y_c^l(t)$ by constraint (29). As a result, the movement of the subject CAVs will affect the upstream traffic flow in the cell in the next time step $f_c^l(t+1)$ and cause congestion.

To model the interaction between the microscopic and macroscopic traffic flow dynamics, we need to recognize the boundary cells that $s$-CAVs, $n$-CAVs and neighbor HDVs belong at each time step. With cell length $\Delta L$, we locate the coordinates for each cell by $\bar{z}_c^l = (\bar{x}_c^l, \bar{y}_c^l)$ and $\underline{z}_c^l = (\underline{x}_c^l, \underline{y}_c^l)$ as Figure 2 shows. Based on these coordinates, (32) and (33) specify the criteria to find the cells involving subject CAVs and neighboring vehicles. Specifically, the binary variable $\phi_{i,c}^l(t)$ takes value one when the coordinates of subject CAV($i$), $z_i(t)$, fall within the cell $c$ in lane $l$ at time $t$, and 0 otherwise. Similarly, $\phi_{\omega,c}^l(t)$ takes value one when the coordinates of neighbor vehicle $\omega$, $\hat{x}_\omega(t)$, fall within the cell $c$ in lane $l$ at time $t$, and 0 otherwise. Moreover, we have constraint (34) to ensure that each vehicle can only be present in one cell at a time. Equation (35) gives the formulation to calculate traffic density $k_c^l(t)$ in a cell. To be noted, it counts both aggregated flow and individual vehicles with trajectory control or estimation in this study. Thus, by combining CTM in (28) - (31) with the cell membership recognition in (32)-(34), (35) can capture the interaction between macroscopic and microscopic traffic dynamics. The logic functions in (32) and (33) can be further transferred to integer constraints given in Appendix A.

$$\phi_{i,c}^l(t) = [z_i(t) \geq \underline{z}_c^l] \wedge [z_i(t) \leq \bar{z}_c^l], \quad \forall\, i \in I, c \in C_l, l \in L \qquad (32)$$

$$\phi_{\omega,c}^l(t) = [\hat{x}_\omega(t) \leq \bar{x}_c^l] \wedge [\hat{x}_\omega(t) \geq \underline{x}_c^l], \quad \forall \omega \in \Omega^l, c \in C_l, l \in L \qquad (33)$$

$$\sum_{l=1}^{L} \sum_{c=1}^{C} \phi_{j,c}^l(t) = 1, \quad \forall j \in J \qquad (34)$$

$$k_c^l(t) = f_c^l(t) + \sum_{j \in J} \phi_{j,c}^l(t) / \delta \Delta L, \quad \forall c \in C_l, l \in L \qquad (35)$$

## 5. MPC for Synchronization Control

MPC is a closed-loop feedback control by repeating the following procedure at each discrete control time step. It takes current upstream traffic's and vehicles' states (both s-CAVs and neighboring HDVs and CAVs) as inputs, predicts their future states in T time steps by vehicle/traffic dynamics models, and then generates the optimal control law for $s$-CAVs in T steps, but only implements the first step, by solving an optimization model. This recursive process will reset the surrounding traffic condition by recognizing new neighbor vehicles and their trajectory profiles at each control sampling step. To a certain extent, it will address the issues regarding prediction error and traffic uncertainty in control. Therefore, it fits this study very well. The MPC built for this study involves a mixed-integer optimizer with multiple objectives to generate the optimal trajectory control law co-considering synchronization and traffic efficiency. By noticing the complexity of weighting control objectives, we further transfer it to two sequential MPCs according to the problem features. Built upon that, we develop an adaptive weighting strategy for balancing the synchronization and traffic efficiency under different states. Below we introduce the technical details for the model development and analysis.

**5.1 Multi-objective MPC and a hybrid MPC system**

The MPC for the synchronization control takes the initial states of $s$-CAVs, $n$-CAVs and the upstream traffic as inputs at time step $t \in \mathbb{Z}^{+}$, and then seeks to generate the optimal control law to instruct subject CAVs' trajectories at $t+1$ so that they can approach each other and synchronize their motion efficiently without jeopardizing traffic efficiency. The synchronization control is subject to the prediction of $s$-CAVs' and $n$-CAVs' movements and the upstream traffic propagation in future $T$ steps according to the micro and macro dynamics introduced in Section 4. The above thoughts are demonstrated by The MPC in (36). Without loss of generality, we make $t$ start with zero in the formulation.

$$J(u_t) = min \sum_{t=0}^{T-1} \frac{1}{2} \left[ \underbrace{u_t^T Q_u u_t}_{J_u} + \underbrace{\mathcal{w}_t^T Q_{\mathcal{w}} \mathcal{w}_t}_{J_{\mathcal{w}}} - \underbrace{\mathcal{y}_t^T Q_{\mathcal{y}} \mathcal{y}_t}_{J_{\mathcal{y}}} - \underbrace{\eta_t^T Q_\eta \eta_t}_{J_\eta} + \underbrace{\mathcal{v}_t{}^T Q_{\mathcal{v}} \mathcal{v}_t}_{J_{\mathcal{v}}} \right] + \underbrace{q_z z^T}_{J_z} \tag{36}$$

For, $t \in T$

s.t. (7) - (10), (12) --(21), and (25) - (35).

Where,

- $u_t = [u_1(t), u_2(t)]^{\mathrm{T}}$ represents the acceleration/deceleration instructions to the two $s$-CAVs at control sample step $t$;
- $\mathcal{w}_t = [\Delta v_1(t), \Delta v_2(t)]^{\mathrm{T}}$ and $\Delta v_i(t)$ represents speed deviations of $s$-CAVs ($i$ =1 or 2) from the maximum speed $\bar{v}$;
- $\mathcal{y}_t$ and $\eta_t$ respectively denote the flow propagation defined by (29) at time step $t$ and the car-following relationship defined by (3).
- $\mathcal{v}$ represents the relative speed of two $s$-CAVs.
- $z = [\Delta x(t) = x_1(t) - x_2(t) - \tilde{d}]_{t=1}^{T}$ represents the deviations of the spacing between the two subject CAVs from the pre-defined desired spacing $\tilde{d}$ during the prediction horizon.
- $Q_u$, $Q_{\mathcal{w}}$, $Q_{\mathcal{y}}$, $Q_\eta$ and $Q_{\mathcal{v}}$ are $2 \times 2$ diagonal weight matrices and $q_z$ is weight vectors, i.e., $Q_{\mathcal{w}} = \frac{1}{(\overline{v_x} - \underline{v_x})^2} \begin{bmatrix} q_{1,\mathcal{w}} & \\ & q_{2,\mathcal{w}} \end{bmatrix}$, which are normalized by the maximum possible speed deviation, $q_{i,\mathcal{w}}$ is a tunable weighting parameter.

The MPC involves six items in the objective function to balance the synchronization and traffic efficiency. Specifically, we evaluate traffic efficiency in this context from three perspectives: (i) the traffic smoothness, (ii) $s$-CAVs' speed, and (iii) upstream traffic propagation. Accordingly, the first item $J_u$ in the objective function puts a penalty on the $s$-CAV's control inputs to promote mild vehicle acceleration/deceleration so that we can keep surrounding traffic smoothness. The second item $J_{\mathcal{w}}$ seeks to sustain $s$-CAVs' speed so that they won't block traffic to facilitate synchronization. The third item $J_{\mathcal{y}}$ maximizes the volume of traffic propagating from the upstream to the downstream cells during the synchronization control. Therefore, the first three items $J_u$ to $J_{\mathcal{y}}$ together help maintain a satisfying traffic efficiency during the synchronization process.

On the other hand, we factor in three aspects to ensure the synchronization efficiency: (i) positioning the $s$-CAV pair in the car-following mode, (ii) stabilizing the spacing of the two $s$-CAVs' to the desired value, and (iii) harmonizing the speed of the two $s$-CAVs. Correspondingly, $J_\eta$ in the objective function rewards the duration that one $s$-CAV is following the other $s$-CAV (i.e., $\eta_{i,i'}^{l}(t) = 1$). The fifth term $J_5$ penalizes on the variation of the relative speed between the two $s$-CAVs. To be noted, $J_{\mathcal{v}}$ benefits both synchronization efficiency and traffic smoothness. And the last item $J_z$ penalizes on the spacing deviation to push the two $s$-CAVs approach and then form a platoon efficiently. Table 1 summarizes the six factors and the corresponding weights.

**Table 1** Objective items and corresponding weights

| **Objective item (weight)** | | **Explanation** |
|---|---|---|
| Traffic performance metric | $J_u$ $(Q_u)$ | Featuring traffic smoothness by $s$-CAVs' control inputs. |
| | $J_w (Q_w)$ | Featuring traffic efficiency by $s$-CAVs' speed. |
| | $J_y$ $(Q_y)$ | Featuring traffic efficiency by upstream flow propagation. |
| Synchronization performance metric | $J_\eta$ $(Q_\eta)$ | Featuring synchronization efficiency by two $s$-CAVs' positions (in car-following mode) |
| | $J_v$ $(Q_v)$ | Featuring speed harmonization by the relative speed of the two $s$-CAVs. |
| | $J_z$ $(q_z)$ | Featuring synchronization efficiency by the relative spacing of $s$-CAVs. |

Clearly, we should design the weight matrices/vector $Q_u$, $Q_w$, $Q_y$, $Q_\eta$, $Q_v$ and $q_z$ carefully to achieve a satisfying control performance. However, this is not trivial. First, the items associated with traffic efficiency and synchronization efficiency are correlated, and to some extent, conflict with each other. For instance, when the two $s$-CAVs are far apart from each other, a good way to expedite the synchronization is to have the leading $s$-CAV slow down and wait for the following $s$-CAV. But it will block traffic and reduce the traffic efficiency. This means if we over prioritize the spacing deviation between the two $s$-CAVs measured in $J_z$, the traffic speed featured by $J_w$ and upstream traffic efficiency featured by $J_y$ will be sacrificed. Moreover, it is also difficult to determine the proper weights for $J_\eta$ to $J_z$ even though they both contribute to synchronization efficiency. Again, consider a scenario in which the two $s$-CAVs are far apart, but drive in different lanes. If we overweigh $J_z$ for reducing the $s$-CAV pair's longitudinal distance, the latitudinal movement will have a risk of being suppressed. As a result, both $s$-CAVs will speed up and approach each other quickly, but still stay in different lanes. In the meantime, we notice that the two $s$-CAVs will physically experience two dynamic traffic states below during the synchronization.

(i) State $q_1$: catch-up state, during which the two s-CAVs on different lanes seek to approach each other and form a car-following mode.

(ii) State $q_2$: platooning state, during which the two s-CAVs try to stabilize their spacing and relative speed to form a platoon on the same lane

Under these two states, the synchronization control has different priorities for the items in the objective function. Under state $q_1$, the $s$-CAVs are far apart and drive in different lanes. It is of greater importance to move these two $s$-CAVs closer to the same lane and position them in car-following mode ($J_\eta$) than harmonizing their speed and spacing ($J_v$ and $J_z$). Under state $q_2$ in which the $s$-CAVs are already under the car-following mode, harmonizing their speed and stabilizing spacing becomes the main focus.

Considering the above problem feature and complexity, we further transfer the multi-objective MPC in (36) into a hybrid MPCs consisting of two sequential MPCs. Specifically, when the pair of $s$-CAVs is initially far away from each other (i.e., under state $q_1$). Accordingly, MPC-I in (37) is used to conduct the synchronization control, only including the objectives in $J_u$ to $J_\eta$ (excludes $J_v$ and $J_z$) to encourage one $s$-CAV catching up the other promptly and then following it in the same lane. Under this stage, stabilizing spacing and relative speed are not our emphases.

**MPC-I**

$$Min\,\Gamma_1(u_t) = \sum_{t=0}^{T-1} \frac{1}{2}\left[\underbrace{u_t^T Q_u u_t}_{J_u} + \underbrace{w_t^T Q_w w_t}_{J_w} - \underbrace{y_t^T Q_y y_t}_{J_y} - \underbrace{\eta_t^T Q_\eta \eta_t}_{J_\eta}\right] \tag{37}$$

For, $t \in T$

s. t. (7) - (10), (12) - (21), and (25) - (35).

Once the following $s$-CAV successfully catches up the other $s$-CAV and they start to run in a car-following mode in the same lane, it triggers the switching signal, by which the two $s$-CAVs switch to the platooning state $q_2$. Accordingly, the synchronization control switches to the MPC-II in (38). It focuses on maintaining the spacing of $s$-CAVs to the desired spacing and synchronize their speed. Accordingly, the optimizer in (38) excludes the term $J_\eta$ in (36) with additional constraint $\sum_{l \in L} \eta_{i,i'}^l(t) = 1$ to ensure $s$-CAVs stay in the same lane.

**MPC-II**

$$Min\,\Gamma_2(u_t) = \sum_{t=0}^{T-1} \frac{1}{2}\left[\underbrace{u_t^T Q_u u_t}_{J_u} + \underbrace{w_t^T Q_w w_t}_{J_w} - \underbrace{y_t^T Q_y y_t}_{J_y} + \underbrace{v_t^T Q_v v_t}_{J_v}\right] + \underbrace{q_z z^T}_{J_z} \tag{38}$$

For, $t \in T$

s. t. $\sum_{l \in L} \eta_{i,i'}^l(t) = 1,\ \forall\, i, i' \in I, t \in T$

(7) - (10), (12) - (21), and (25) - (35).

Note that both MPC-I and MPC-II want to sustain traffic efficiency, even though they emphasize on different aspects of synchronization efficiency. Therefore, both of them involve speed fluctuation ($J_u$), speed deviation ($J_w$) and upstream traffic efficiency ($J_y$) in the objective functions. The feasibility of the hybrid MPC system is analyzed in Section 6.

### 5.2 Weighting Strategies

Even though the hybrid MPC system reduces the weighting complexity in the objective function by decorrelating the items in (36) regarding the synchronization efficiency into two MPCs, each is still a multi-objective optimizer, aiming to balance the traffic and synchronization efficiency in the respective state. Moreover, we noticed that control priority may vary according to the change of the initial traffic state at each control time stamp. Taking the catch-up state ($q_1$) as an example, we explain this thought by the two initial states: (i) two $s$-CAVs drive with normal speed in two non-adjacent lanes; (ii) two $s$-CAVs drive with slow speed in two adjacent lanes. Recall that $J_w$ in the objective function of MPC-I in (37) aims to sustain $s$-CAVs' speeds without causing adverse traffic impact, and $J_\eta$ seeks to position $s$-CAVs in the car-following mode. Clearly, each state favors either $J_w$ or $J_\eta$ and it encourages us to use different weight schemes toward $J_w$ and $J_\eta$ for the best synchronization control. More exactly, we want to award $J_\eta$ more under initial state (i) to position the two $s$-CAVs in the car-following mode quickly but give priority to $J_w$ under the initial state (ii) to speed up s-CAVs and mitigate their adverse traffic impact. One the other hand, we found that only the weights of $J_w$ and $J_\eta$ in MPC-I and the weights of $J_w$ and $J_z$ in MPC-II need to be adaptively adjusted according to the initial states, while $J_u$, $J_y$ (and $J_v$), respectively facilitating traffic smoothness, upstream traffic efficiency, (and $s$-CAV speed harmonization), should be consistently considered under any initial states when the control is in the catch-up state ($q_1$) (and platooning state ($q_2$)). Therefore, the adaptive weighting strategies will mainly focus on the tradeoff between $J_w$ and $J_\eta$ in MPC-I ( $J_w$ and $J_z$ in MPC-II).

To develop the weighting strategies for MPC-I, we first introduce the optimality loss regarding $J_w$ and $J_\eta$. It should be noted that $J_w$ in the objective function (37) penalizes $s$-CAVs' speed deviation from the maximum speed, $\bar{v}_x$. Thus, it reaches to the minimum value when the $s$-CAVs drive at the maximum speed. In the meantime, $J_\eta$ in the objective function (37) rewards the duration that the two $s$-CAVs are running under the car-following mode, and it gains more values when $s$-CAVs drive in the same lane. Invoked by these observations, we first introduce (39) to measure the optimality loss of $J_w$ under a given control law.

$$\Delta J_{w,i}(t) = \frac{\bar{v}_x t\delta - x_i(t)}{\bar{v}_x t\delta}, t \in \mathbb{Z}^+ \tag{39}$$

where $\bar{v}_x t\delta$ represents the furthest longitudinal position that $s$-CAV$(i)$ can reach up to the time step $t$ if it drives at the maximum speed $\bar{v}_x$; $x_i(t)$ is the actual position of $s$-CAV$(i)$ at step $t$. Then, $\Delta J_{w,i}(t)$ represents the normalized loss of $s$-CAV's longitudinal displacement and it characterizes the optimality loss of the control regarding the objective in $J_w$. A smaller $\Delta J_{w,i}(t)$ indicates a smaller optimality loss of $J_w$ and thus, a smaller weight can be put on $J_w$, and vice versa.

Following a similar idea, we use (40) to characterize the optimality loss of the reward in $J_\eta$ under a given control law.

$$\Delta J_\eta(t) = \frac{|y_1(t) - y_2(t)|}{\bar{y}^L - \underline{y}^L}, t \in \mathbb{Z}^+ \tag{40}$$

where $|y_1(t) - y_2(t)|$ is the lateral spacing between the two $s$-CAVs and $\bar{y}^L - \underline{y}^L$ is the road width. Clearly, the control law gains small rewards from $J_\eta$ in the objective funcaiton (i.e., $\Delta J_\eta$ is large) when $s$-CAVs move in different lanes, but a large and accumulated reward when the two $s$-CAVs run in the same lane (i.e., $\Delta J_\eta = 0$).

Using the two optimality losses, we develop the weighting strategy for MPC-I in (41), which adaptively tunes the weights of $J_w$ and $J_\eta$ in (37) aiming to have the following $s$-CAV to quickly catch up the leading $s$-CAV and then follow the car-following mode in the same lane, while sustaining traffic efficiency.

$$q_\eta = \sum_{i \in I} \alpha_i \, q_{i,w} \xi_{i,t}^I, \; \xi_{i,t}^I = \frac{\Delta J_\eta}{\Delta J_{w,i}} \tag{41}$$

where $q_{i,w}$ and $q_\eta$ are the weight assigned to objective $J_w$ and $J_\eta$, respectively. $\alpha_i$ is a tunable parameter, and $\xi_{i,t}^I$ is a scaling parameter equal to the ratio of the optimality loss of $J_\eta$ to that of $J_w$. The key idea in (41) is that if the two $s$-CAVs are moving in different lanes (i.e., a large $\Delta J_\eta(t)$), we give a higher priority to $J_\eta$ in (37) through varying the scaling parameter $\xi_{i,t}^I$ to encourage lane changes which moves two $s$-CAVs in the same lane and position them in a car-following mode unless they drive too low (i.e., a large $\Delta J_{w,i}(t)$), i.e., weigh $J_\eta$ over $J_w$ in the objective function, and vice versa.

To develop the weighting strategy for MPC-II, we notice that the last item $J_z$ in the objective function (38) focus on stabilizing the spacing between two $s$-CAVs (i.e., penalizing the deviation of $s$-CAVs' spacing from the desired spacing). Following the same idea as before, we introduce (42) to measure the optimality loss regarding the objective $J_z$ under a given control law.

$$\Delta J_z(t) = \frac{|x_2(t) - x_1(t) - \tilde{d}|}{\tilde{d}} \tag{42}$$

where $|x_2(t) - x_1(t) - \tilde{d}|$ is the deviation of $s$-CAVs' spacing from the desired spacing, $\tilde{d}$. Then, $\Delta J_z(t)$ represents the normalized spacing deviation of $s$-CAVs and it characterizes the optimality loss of the control regarding the objective $J_z$. A larger $\Delta J_z(t)$ indicates a larger optimality loss of $J_z$ (larger deviation from desired spacing) and thus, $J_z$ can be prioritized, and vice versa.

Using the two optimality losses $\Delta J_{w,i}$ in (39) and $\Delta J_z$ in (42), we develop the adaptive weighting strategy for MPC-II in (43), aiming to keep the two $s$-CAVs in an efficient platooning mode, while sustaining the traffic efficiency.

$$q_z = \sum_{i \in I} \alpha_i \, q_{i,w} \xi_{i,t}^{II}, \xi_{i,t}^{II} = \frac{\Delta J_z}{\Delta J_{w,i}} \tag{43}$$

where $q_{i,w}$ and $q_z$ are the weights assigned to objective $J_w$ and $J_z$, respectively. $\alpha_i$ is a tunable parameter and $\xi_{i,t}^{II}$ is a scaling parameter equal to the ratio of the optimality loss of $J_z$ to that of $J_w$. Here, the main idea is that if the spacing of the two s-CAVs deviates too much from the desired spacing (a large $\Delta J_z$), we give a high priority to $J_z$ in the objective function (38) through the scaling parameter $\xi_{i,t}^{II}$, which encourages $s$-CAVs to adjust their speed and achieve the desired spacing unless they drive too slow and block the traffic (a large $\Delta J_{w,i}$), and vice versa.

We also notice that $\Delta J_{w,i} = 0$ is a feasible scenario for both MPC-I and II, and it will raise the numerical issues to apply (42) and (43). To address this problem, we limit the scaling parameters $\xi_{i,t}^{I}$ (and $\xi_{i,t}^{II}$) to a numerical range, i.e., $\xi_{i,t}^{I} \in [0, \xi_{max}^{I}]$, where $\xi_{max}^{I}$ is a given parameter to bound $\xi_{i,t}^{I}$ and avoid numerical issue. In addition, the scaling parameter $\xi_{i,t}^{I}$ (and $\xi_{i,t}^{II}$) will be calculated at each sample time point $t$ and keep constant within a prediction horizon.

## 6. Recursive Feasibility of the Constraint Set

As the MPC is implemented recursively at each time step, a fundamental question is whether the MPC can find a feasible control law at each time step $t$ (i.e., whether the constraint set of the MPC optimizer is non-empty at each time step $t$), given the platoon system starts from an initial feasible condition. A system is called recursive (persistent) feasible (Mayne, 2014) if the answer to this question is affirmative. Apparently, it is extremely important for the applicability of the MPC in practical implementation. Therefore, this section will investigate the recursive feasibility of the MPC-I (37) and MPC-II (38). Before that, we need to clarify two points. Recall the MPC will reset the initial condition at each control time step by recognizing new neighbor vehicles $\Omega$, their positions and velocities. Therefore, the MPC accounts for the uncertainty and dynamics of neighbor traffic and surrounding vehicles. A new control law at each step can correct control errors resulting from previous steps due to inaccurate predictions. Moreover, given MPC-I and MPC-II only involve two s-CAVs and a limited number of neighbor vehicles, their model sizes are relatively small and can be efficiently solved by a commercial solver, i.e., Gurobi, within the control interval. Therefore, the control continuity is sustained. These views are fundamentals and prerequisites of the recursive feasibility analysis.

Given MPC-I and MPC-II share the same constraint set, we will only prove that the recursive feasibility of MPC-I. To do that, we first formally define this constraint set $\mathcal{P}(\boldsymbol{u}(t), \boldsymbol{\chi}(t))$ by (44), and want to prove that if $\mathcal{P}(\boldsymbol{u}(t-1), \boldsymbol{\chi}(t-1))$ is nonempty at control step $t-1$, then the constraint set $\mathcal{P}(\boldsymbol{u}(t), \boldsymbol{\beta}(t))$ is nonempty at control step $t$, $t \in \mathbb{Z}_+$.

$$\mathcal{P}(\boldsymbol{u}(t), \boldsymbol{\chi}(t)) = \{\boldsymbol{u}(t) \in \mathbb{R}^{T|I|}, \boldsymbol{\chi}(t) \in \mathbb{B}^{nT} \,|\, g_t(\boldsymbol{u}(t), \boldsymbol{\chi}(t)) \leq 0 \}, \forall t \in \mathbb{Z}_+ \tag{44}$$

where $\mathbb{B} = \{0,1\}$, and $g_t(\boldsymbol{u}(t), \boldsymbol{\chi}(t)) \leq 0$ is the constraint set defined by constraints (7) - (10), (12) - (21), and (26) - (35). $\boldsymbol{\chi}(t)$ is the set of binary variables, respectively including vehicles' lane memberships, relative positions of vehicle pairs and cell memberships, i.e., $\boldsymbol{\chi}(t) = \{\boldsymbol{\gamma}(t), \boldsymbol{\rho}(t), \boldsymbol{\eta}(t), \boldsymbol{\phi}(t)\}$.

The standard approach to prove the recursive feasibility of an MPC is to show there exists a control-invariant set that contains the origin. However, it is generally difficult to prove the recursive feasibility of an MPC with binary variables using the standard approach. This is because the combinations of binary variables lead to an enormous number of solution spaces, and it is hard to show the existence of a control-invariant set among them (Marcucci and Tedrake, 2020). However, the constraint set $\mathcal{P}(\boldsymbol{u}(t), \boldsymbol{\chi}(t))$ in this study presents unique features as follows which can facilitate our proof. (i) The values of binary variables in $\boldsymbol{\chi}(t)$ are determined by the leading (following) relationships of the two $s$-CAVs and their neighbor vehicles, which presents a limited number of cases (will define them in detail in the proof later). Accordingly, the constraint set (44) can be divided into small constraint sets that only contain the continuous variables $\boldsymbol{u}(t)$. (ii) Constraints (28) - (35), which describe the upstream traffic propagation, are naturally feasible once vehicles' cell membership, i.e., binary variables $\boldsymbol{\phi}(t)$, are determined. Thus, we only consider constraints (7) - (10), (12) - (21), (26), and (27), and omit constraints (28) - (35) in our proof. Taking advantage of these features, we prove the recursive feasibility of MPC-I in **Theorem 1**.

Mainly, **Theorem 1** claims that MPC-I is feasible at step $t$, i.e., $\mathcal{P}(\boldsymbol{u}(t), \boldsymbol{\beta}(t)) \neq \emptyset$ if the four conditions are satisfied. Specifically, Condition (i) calls for the feasibility of MPC-I at step $t-1$, i.e., $\mathcal{P}(\boldsymbol{u}(t-1), \boldsymbol{\chi}(t-1)) \neq \emptyset$. Condition (ii) requires that neighbor vehicles follow the speed limits as the subject CAVs, i.e., $\hat{v}_\omega(t) \in \left[\underline{v_x}, \overline{v_x}\right], \forall \omega \in \Omega_h \cup \Omega_c$. This is a reasonable condition since in a normal traffic, every vehicle on the road follows the regulatory speed limits. Condition (iii) requires the reaction time of neighbor HDVs to satisfy $\tau_\omega \leq \bar{\tau}_\omega, \forall \omega \in \Omega_h$. We will demonstrate the rationality of the bound in the remark followed by our proof of **Theorem 1**. Condition (iv) requires that when an $n$-CAV($\omega$) follows a $s$-CAV($i$), the relative speed between them satisfies $v_i(t) - \hat{v}_\omega(t) \geq \frac{\tilde{a}_i(t-1)\delta\left(B-\frac{1}{2}\right)}{A}, \forall\, \omega \in \Omega_c$. In other words, the speed of $s$-CAV, $v_i(t), \forall i \in I$, is bounded by $v_i(t) \geq \hat{v}_\omega(t) + \frac{\tilde{a}_i(t-1)\delta\left(B-\frac{1}{2}\right)}{A}$. We will show that this condition is recursively satisfied if it holds at the beginning of the prediction horizon, i.e., $t = 0$, in **Theorem 1**. Moreover, by plugging in the values of parameters (Xiao et al., 2017), i.e., $\delta = 1\, s, A = 0.1836, B = 0.8163, \tilde{a}_i(t-1) = \underline{a_x} = -6\, m/s^2$, we can obtain that $\frac{\tilde{a}_i(t-1)\delta\left(B-\frac{1}{2}\right)}{A} = -10.33\, m/s$ (or $-23.11\, mph$). It indicates that Condition (iv) is satisfied in most of the normal traffic scenarios.

**Theorem 1.** $\mathcal{P}(\boldsymbol{u}(t), \boldsymbol{\chi}(t))$ is recursively feasible, i.e., $\mathcal{P}(\boldsymbol{u}(t), \boldsymbol{\chi}(t)) \neq \emptyset, \forall t \in \mathbb{Z}_+$, if the following conditions hold.

(i) $\mathcal{P}(\boldsymbol{u}(t-1), \boldsymbol{\chi}(t-1)) \neq \emptyset, \forall t \in \mathbb{Z}_+$.

(ii) $\underline{v_x} \leq \hat{v}_\omega(t) \leq \overline{v_x}, \forall t \in \mathbb{Z}_+, \omega \in \Omega_h \cup \Omega_c$.

(iii) $\tau_\omega \leq \bar{\tau}_\omega, \forall \omega \in \Omega_h$, *where* $\bar{\tau}_\omega = \frac{-(\underline{v_x}-5\underline{a_x})-\sqrt{(\underline{v_x}-5\underline{a_x})^2-24\underline{a_x}d_\omega}}{2\underline{a_x}}, \omega \in \Omega_h$

(iv) *The speeds of $n$-CAVs($\omega$) and its leading vehicle, $i$, satisfy* $v_i(t) - \hat{v}_\omega(t) \geq \frac{\tilde{a}_i(t-1)\delta(B-\frac{1}{2})}{A}$, *where* $\tilde{a}_i(t) = max\left[\underline{a_x}, \underline{a_{i,v}}(t)\right]$ *and* $\underline{a_{i,v}}(t) = \frac{\underline{v_x}-v_i^x(t)}{\delta}$.

**Proof.** To begin with, we have two $s$-CAVs in the study and we noticed that except for the safety constraint (12), in which two $s$-CAVs are coupled (will discuss later in detail), the constraints in $\mathcal{P}(\boldsymbol{u}(t), \boldsymbol{\chi}(t))$ are formulated for each $s$-CAV separately. Thus, we can separate the constraint set such that $\mathcal{P}(\boldsymbol{u}(t), \boldsymbol{\chi}(t)) = \{\mathcal{P}(u_1(t), \boldsymbol{\chi}(t)), \mathcal{P}(u_2(t), \boldsymbol{\chi}(t))\}$, where $u_1(t)$ and $u_2(t)$ are respectively the control inputs for $s$-CAV(1) and $s$-CAV(2). Without loss of generality, our proof below will first focus on one $s$-CAV such as $s$-CAV(1) and then the results can be easily extended to the other $s$-CAV so that we can conclude **Theorem 1.**

To conduct the proof, we further rewrite the constraint sets $\mathcal{P}(u_1(t), \boldsymbol{\chi}(t))$ in (45) according to the binary variables within the constraints.

$$\mathcal{P}(u_1(t), \boldsymbol{\chi}(t)) = \mathcal{P}_u(u_1(t)) \cap \mathcal{P}_\chi(u_1(t), \boldsymbol{\chi}(t)) \tag{45}$$

Where,

(1) $\mathcal{P}_u(u_1(t))$: The constraint set only includes continuous variables in Equation (7) - (10) and Equation (17) - (19) for capturing vehicle dynamics of $s$-CAV(1) at step $t \in T$.
(2) $\mathcal{P}_\chi(u_1(t), \boldsymbol{\chi}(t))$: The constraints set includes both binary variables $\boldsymbol{\chi}(t)$ and continuous variables in constraints (12), (13), (15), (16), (20), (21), (26) and (27) for capturing the following relationship and safety constraints of $s$-CAV(1) at step $t \in T$.

Next, we discuss the three cases for specifying the following relationships between $s$-CAV(1) and its leading or following neighbor vehicles, each of which determines a set of values of the binary variables in $\boldsymbol{\chi}(t)$. Case (1) considers that $s$-CAV(1) follows another vehicle $j$ in lane $l^* \in L$, which can be the other $s$-CAV, or a $n$-CAV, or an HDV. A vehicle pair in this case is denoted as $1 \to j, \forall j \in J$. According to Equations (1) - (6), Case (1) determines the values of the lane membership indicators ($\gamma$) and the following relationship indicators ($\eta$ and $\rho$). We use $\chi_1(u_1(t))$ to denote this set of solutions shown in Equation (46). Then, the constraints set $\mathcal{P}_\chi(\boldsymbol{u}(t), \boldsymbol{\chi}(t))$ is simplified into the constraints associated with the vehicle pair $1 \to j$, which we denote as $\mathcal{P}_\chi(u_1(t)|\chi_1)$.

$$\begin{aligned}\chi_1(u_1(t)) = \{&\gamma_1^{l^*}(t) = 1, \gamma_j^{l^*}(t) = 1, \eta_{1,j}^{l^*}(t) = 1, \rho_{1,j}^{l^*}(t) = 1\} \\ &\cup \{\gamma_1^{l}(t) = 0, \gamma_j^{l}(t) = 0, \eta_{1,j}^{l}(t) = 0, \rho_{1,j}^{l}(t) = 0\}_{l \in L \backslash l^*}, \forall j \in J, l^* \in L\end{aligned} \tag{46}$$

For Case (2), we consider that $s$-CAV(1) is followed by an HDV($\omega$) in lane $l^* \in L$. We denote this vehicle pair as $\omega \to 1, \forall\, \omega \in \Omega_h$. Following the similar discussion used in Case (1), we determine the values of binary variables as shown in (47) according to Equations (1) - (6), and then simplify the constraints set $\mathcal{P}_\chi(\boldsymbol{u}(t), \boldsymbol{\chi}(t))$ to $\mathcal{P}_\chi(u_1(t)|\chi_2)$.

$$\begin{aligned}\chi_2(u_1(t)) = \{&\gamma_\omega^{l^*}(t) = 1, \gamma_1^{l^*}(t) = 1, \eta_{\omega,1}^{l^*}(t) = 1, \rho_{\omega,1}^{l^*}(t) = 1\} \\ &\cup \{\gamma_1^{l}(t) = 0, \gamma_\omega^{l}(t) = 0, \eta_{1,\omega}^{l}(t) = 0, \rho_{1,\omega}^{l}(t) = 0\}_{l \in L \backslash l^*}, \forall\, \omega \in \Omega_h, l^* \in L\end{aligned} \tag{47}$$

Similarly, in Case (3), we consider that $s$-CAV(1) is followed by a $n$-CAV($\omega$) in lane $l^* \in L$ and denote this vehicle pair as $\omega \to 1, \forall \omega \in \Omega_c$ . Then, the binary variables' values in shown (48) can be determined by Equations (1) - (6). We denote $\mathcal{P}_\chi(u_1(t)|\chi_3)$ as the activated constraints associated with these binary variables in $\mathcal{P}_2$.

$$\begin{aligned}\chi_3(u_1(t)) = \{&\gamma_\omega^{l^*}(t) = 1, \gamma_1^{l^*}(t) = 1, \eta_{\omega,1}^{l^*}(t) = 1, \rho_{\omega,1}^{l^*}(t) = 1\} \\ &\cup \{\gamma_1^{l}(t) = 0, \gamma_\omega^{l}(t) = 0, \eta_{1,\omega}^{l}(t) = 0, \rho_{1,\omega}^{l}(t) = 0\}_{l \in L \backslash l^*}, \forall \omega \in \Omega_c, l^* \in L\end{aligned} \tag{48}$$

The three cases specify the leading and following relationship of s-CAV(1) and its neighbors at any time step $t$. Specifically, it will always drive under Case (1) combined with either Case (2) or Case (3).

Then, we reformulate the constraint set $\mathcal{P}_\chi$ for $s$-CAV(1) in (49). For example, if $s$-CAV(1) follows a $n$-CAV, and is followed by an HDV, then $\mathcal{P}_\chi$ is specified as $\mathcal{P}_\chi(u_1(t)|\chi_1) \cup \mathcal{P}_\chi(u_1(t)|\chi_2)$.

$$\mathcal{P}_\chi\big(u_1(t), \boldsymbol{\chi}(t)\big) = \mathcal{P}_\chi(u_1(t)|\chi_1) \cup \mathcal{P}_\chi(u_1(t)|\chi_2) \cup \mathcal{P}_\chi(u_1(t)|\chi_3) \tag{49}$$

Built upon (49), we want to show $\mathcal{P}_u(u_1(t)) \cap \mathcal{P}_\chi(u_1(t)|\chi_1) = \mathbb{U} \neq \emptyset$, and then demonstrate $\mathbb{U} \cap \mathcal{P}_\chi(u_1(t)|\chi_2) \neq \emptyset$ and $\mathbb{U} \cap \mathcal{P}_\chi(u_1(t)|\chi_3) \neq \emptyset$, given the four conditions in Theorem 1 are satisfied. This proof can be extended to $s$ $-$CAV(2) so that we can draw the conclusion and complete the proof for Theorem 1. Below are the derivation details for each step.

**Scenario I**: without loss of generality, we start with the scenario that $s$-CAV(1) follows an HDV($\omega$) in a lane $l^* \in L$ (Case(1)). By plugging the values of the integer variables given in (46) into Constraints (12), (13), (15), (16), (20), (21), (26), and (27), we obtain the only activated constraint (13) in $\mathcal{P}_\chi(u_1(t)|\chi_1)$, for which we introduce a new notation $\hat{\wp}_{1,\omega}\big(u_1(t)\big)$ as follows.

$$\hat{\wp}_{1,\omega}\big(u_1(t)\big) := \hat{x}_\varpi(t) - x_1(t) - (L_1 + \tau_1 v_1^x(t) - \frac{\left[v_1^x(t) - \underline{v_x}\right]^2}{2\underline{a_x}}), \forall\ \omega \in \Omega_h \tag{50}$$

Then, combining constraints (17) and (18) in $\mathcal{P}_u\big(u_1(t)\big)$, we ensure that the speed of HDV($\omega$) follows $\underline{v_x} \leq \hat{v}_\omega(t) \leq \overline{v_x}$, which confirms Condition (ii) in **Theorem 1**. Now, $\mathcal{P}_u\big(u_1(t)\big) \cap \mathcal{P}_\chi(u_1(t)|\chi_1) = \mathbb{U}(u_1(t))$ is the joint constraint set involving constraints (7) - (10) and safety constraint (50). We claim that $\mathbb{U}(u_1(t))$ is non-empty and given in (51) for any $\hat{v}_\omega(t) \in [\underline{v_x}, \overline{v_x}]$ by referring to the Lemma 4.1 in (Gong et al., 2016) without providing the detailed proof.

$$u_1(t) \in \mathbb{U}\big(u_1(t)\big) = [max\left\{\underline{a_x}, \underline{a_{1,v}}(t)\right\}, min\{\overline{a_x}, \overline{a_{1,v}}(t), \overline{a_{1,d}}(t)\}] \tag{51}$$

Where,

$$\underline{a_{1,v}}(t) = \frac{\underline{v_x} - v_1^x(t)}{\delta} \leq 0 \tag{51.1}$$

$$\overline{a_{1,v}}(t) = \frac{\overline{v_x} - v_1^x(t)}{\delta} \geq 0 \tag{51.2}$$

$$\overline{a_{1,d}}(t) = \frac{3}{2}\underline{a_x} + \underline{a_{1,v}} - \frac{\underline{a_x}}{\delta^2}\sqrt{B(t)} \tag{51.3}$$

$$B(t) = \frac{\left(v_1^x(t-1) - \underline{v_x}\right)^2 \delta^2}{\underline{a_x}^2} + \frac{\delta^3 (v_1^x(t-1) - \underline{v_x})}{-\underline{a_x}} + \frac{9}{4}\delta^4 + \left(\frac{2\delta^2}{-\underline{a_x}}\right)\left[\left(\frac{\hat{v}_\omega(t-1) + \hat{v}_\omega(t)}{2} - \underline{v_x}\right)\delta + \hat{\wp}_{1,\omega}\left(u_1(t-1)\right)\right] \geq 0 \tag{51.4}$$

Built upon $\mathbb{U}(u_1(t))$, we next consider that this $s$-CAV(1) is followed by an HDV($\omega$), $\forall \omega \in \Omega_h$ (i.e., a scenario formed by Case (1) combined with Case (2)). By plugging the values of the binary variables given in (47) into constraints (12), (13), (15), (16), (20), (21), (26), and (27), we find that the only activated constraints in $\mathcal{P}_\chi(u_1(t)|\chi_2)$ are Equations (15), (16), which determines the spacing of the HDV. Therefore, the activated constraints in $\mathcal{P}_u(u_1(t)) \cap \mathcal{P}_\chi(u_1(t)|\chi_2)$ are given in the set of (52), where $k = \frac{\tau_\omega}{\delta} \in \mathbb{Z}_+$. Clearly, (52) will not be empty as long as there exists spacings $s_{\omega,1} \geq 0$. Note that $k$ is the number of time steps corresponding to the HDV's response delay to its leading vehicle. When $k \leq 1$ for which the HDV can quickly respond to the control of $s$-CAV(1), we only need to show that there exists a positive spacing $s_{\omega,1}$ at step $t+1$ in (52). When $k > 1$, for which there exists a response delay of HDV that is greater than the step size, we have to examine the constraints in extra time steps, $t+k$ and afterwards. Without loss of generality, we assume $k > 1$, given that the feasibility can be easily shown for $k \leq 1$ if the feasibility holds for $k > 1$.

$$v_1(t+k) = v_1(t) + \delta k u_1(t) \tag{52.1}$$

$$x_1(t+k) = x_1(t) + \delta k v_1(t) + \frac{(\delta k)^2}{2} u_1(t) \tag{52.2}$$

$$\hat{v}_\omega(t+k) = max\left\{\underline{v_x}, \frac{\delta}{\tau_\omega}\left(s_{\omega,1}(t) - d_\omega\right)\right\} \tag{52.3}$$

$$\hat{v}_\omega(t+k) = min\left\{\overline{v_x}, \frac{\delta}{\tau_\omega}\left(s_{\omega,1}(t) - d_\omega\right)\right\} \tag{52.4}$$

$$\hat{x}_\omega(t+k) = \hat{x}_\omega(t) + \tau_\omega \hat{v}_\omega(t) \tag{52.5}$$

$$s_{\omega,1}(t+k) = x_1(t+k) - \hat{x}_\omega(t+k) \tag{52.6}$$

As a result, we seek to prove that there exists $u_1(t) \in \mathbb{U}(u_1(t))$ such that the spacing $s_{\omega,1} \geq 0$ in step $t+k$ and afterward. Given $\mathcal{P}_\chi(u_1(t-1)|\chi_2)$ is feasible, we have $s_{\omega,1}(t-1) \geq 0$. Then we discuss it in two situations as follows.

(i) When $0 \leq s_{\omega,1}(t-1) < \underline{v_x}\frac{\tau_\omega}{\delta} + d_\omega$, according to Equation (52.4), we have $\hat{v}_\omega(t-1+k) = \underline{v_x}$. Given $\hat{v}_\omega(t-1+k) = \underline{v_x} \leq v_1(t-1+k)$, we can expect that $s_{\omega,1}(t-1+k) > s_{\omega,1}(t-1) \geq 0$.

(ii) When $s_{\omega,1}(t-1) \geq \underline{v_x}\frac{\tau_\omega}{\delta} + d_\omega$, we further consider two more subcases: (1) $\underline{v_x}\frac{\tau_\omega}{\delta} + d_\omega \leq s_{\omega,1}(t-1) < \overline{v_x}\frac{\tau_\omega}{\delta} + d_\omega$. It is equivalent to $\underline{v_x} \leq \frac{\delta}{\tau_\omega}\left(s_{\omega,1}(t-1) - d_\omega\right) \leq \overline{v_x}$. Based on (52.3), we have $\hat{v}_\omega(t-1+k) = max\left\{\underline{v_x}, \frac{\delta}{\tau_\omega}\left(s_{\omega,1}(t-1) - d_\omega\right)\right\} = \frac{\delta}{\tau_\omega}\left(s_{\omega,1}(t-1) - d_\omega\right)$, which further indicates $s_{\omega,1}(t-1) = \hat{v}_\omega(t-1+k)\frac{\tau_\omega}{\delta} + d_\omega$. Then, we consider the other subcase (2) when $s_{\omega,1}(t-1) \geq \overline{v_x}\frac{\tau_\omega}{\delta} + d_\omega$. It means $\overline{v_x} \leq \frac{\delta}{\tau_\omega}\left(s_{\omega,1}(t-1) - d_\omega\right)$. According to (52.4), $\hat{v}_\omega(t-1+k) = min\left\{\overline{v_x}, \frac{\delta}{\tau_\omega}\left(s_{\omega,1}(t-1) - d_\omega\right)\right\} = \overline{v_x}$. Therefore, we have $s_{\omega,1}(t-1) \geq \hat{v}_\omega(t-1+k)\frac{\tau_\omega}{\delta} + d_\omega$. Combining the results from the two subcases, we conclude that $s_{\omega,1}(t-1)$ satisfy $s_{\omega,1}(t-1) \geq \hat{v}_\omega(t-1+k)\frac{\tau_\omega}{\delta} + d_\omega$.

Now, based on (52.1) - (52.6), we rewrite the spacing $s_{\omega,i}$ in a recursive form as follows.

$$\begin{aligned} s_{\omega,1}(t+nk) = s_{\omega,1}(t+(n-1)k) - s_{\omega,1}(t+(n-2)k) + \frac{1}{2}u_1(t+(n-1)k)\delta^2 \\ + v_1(t+(n-1)k)\delta + d_\omega, \forall n \in \mathbb{Z}_+ \end{aligned} \tag{53}$$

According to (52.1), when $u_1(t)$ takes the minimum value, $v_1(t+k)$ reaches the minimum value. Then, upon (53), we have $s_{\omega,1}(t+nk)$ reach the minimum value if $u_1(t+n'k)$ takes the minimum value for $\forall n' = 0,1,2\dots,n$. Therefore, in order to show that $s_{\omega,1}(t+nk) > 0, \forall n \in \mathbb{Z}_+$, we consider an extreme scenario where $u_1(t+n'k)$ takes the minimum value (i.e., the maximum deceleration) defined by $\mathbb{U}, \forall n' = 0,1,2\dots,n$. From (51), the potential maximum deceleration in $\mathbb{U}$ is $\underline{a_x}$ and the corresponding smallest spacing defined in (53) is given by $s_{\omega,1}(t+5k) = \hat{v}_\omega(t-1+k)\tau_\omega + d_\omega + 6\underline{a_x}{\tau_\omega}^2$ (see Appendix B for the derivation details). Given that $\tau_\omega \leq \frac{-(\underline{v_x}-5\underline{a_x})-\sqrt{(\underline{v_x}-5\underline{a_x})^2-24\underline{a_x}d_\omega}}{2\underline{a_x}}$ from Condition (iii) of **Theorem 1**, we have $s_{\omega,1}(t+5k) \geq 0$, which confirms that the spacing $s_{\omega,1}(t+nk)$ will be positive for any input $u_1(t) \in \mathbb{U}$, if Condition (iii) holds.

Here we want to consider a state transition and its feasibility further. Namely, $s$-CAV(2) cuts in between $s$-CAV(1) and its follower, HDV($\omega$) at step $t$. Then, the scenario transfers to the special one where $s$-CAV(2) follows $s$-CAV(1) and is followed by the HDV($\omega$). This transition makes the constraint of the two s-CAVs correlated in the safety constraint. Then, we need to prove the recursive feasibility regarding s-CAV(2) in this special scenario. Mainly, we want to show $\mathcal{P}_u(u_2(t)) \cap \mathcal{P}_\chi(u_2(t)|\chi_1) = \mathbb{U}\big(u_2(t)\big) \neq \emptyset$, and then there exists $u_2(t) \in \mathbb{U}\big(u_2(t)\big)$ such that $\mathcal{P}_u(u_2(t)) \cap \mathcal{P}_\chi(u_2(t)|\chi_2) \neq \emptyset$, given Conditions (i)-(iv) in **Theorem 1** hold. The proof is similar to the argument above, thus we only present the main idea as follows. First of all, with the given following relationship, we can specify the relevant integer variables in(46), i.e., $\gamma_1^{l^*}(t)=1, \gamma_2^{l^*}(t)=1, \eta_{2,1}^{l^*}(t)=1, \rho_{2,1}^{l^*}(t)=1$, and then the only activated constraint in $\mathcal{P}_\chi(u_2(t)|\chi_1)$ is constraint (12), which couples two $s$-CAVs. According to constraint (10), $v_2(t) \in [\underline{v_x}, \overline{v_x}]$, we confirm Condition (ii) of **Theorem 1**. Then, following the same argument for $s$-CAV(1) under Case (1), we can also show that $\mathbb{U}\big(u_2(t)\big) \neq \emptyset$. Next, given the cutting-in occurs at time step $t$, we have $\eta_{\omega,2}^{l^*}(t)=1$ and $\eta_{\omega,2}^{l^*}(t-1)=0$ (i.e., HDV($\omega$) starts to follow $s$-CAV(2) at step $t$ but not $t-1$). Accordingly, the safety constraint (20) is "turned on" to ensure the spacing between $s$-CAV(2) and the HDV $\omega$ satisfying $s_{\omega,2}(t) = x_2(t) - \hat{x}_\omega(t) \geq \frac{\tau_\omega}{\delta}\hat{v}_\omega(t) + d_\omega$. Then, through the similar discussion for $s$-CAV(1) under the scenario of Case (1) combined with Case (2), we can prove that $s_{\omega,2}(t+nk) \geq 0$ for any control input $u_2(t) \in \mathbb{U}\big(u_2(t)\big)$. The above arguments together conclude the feasibility after this state transition under Scenario I.

**Scenario II** is the combination of Case (1) and Case (3). Namely, we consider $s$-CAV(1) follow a neighbor vehicle $j \in J$ and it is followed by a $n$-CAV($\omega$) in lane $l^* \in L$. We want to show that there exists $u_1(t) \in \mathbb{U}\big(u_1(t)\big)$ such that $\mathcal{P}_u\big(u_1(t)\big) \cap \mathcal{P}_\chi(u_1(t)|\chi_3) \neq \emptyset$, given that $\mathcal{P}(u_1(t-1), \chi)$ and other conditions in **Theorem 1** hold. As $n$-CAV($\omega$) follows the $s$-CAV(1) at step $t-1$ , we have $\eta_{\omega,1}^{l^*}(t-1) = 1$. Accordingly, Constraint (27) is relaxed. Then the only activated constraint in $\mathcal{P}_\chi(u_1(t)|\chi_3)$ is Equation (26). Given that the CACC scheme is operated with $k_1 = 0.01\ s^{-2}$ , $k_2$= 1.6 $s^{-1}$ and $t_d = 0.6s$ (Xiao et al., 2017). Accordingly, we have $k_2 \gg k_1$ , $A = 0.1836$ and $B = 0.8163$ , and $C \approx 0$ from the formulations of $A$ and $B$ in (25). Therefore, we can drop the items related with $C$ in Equation (26), and then it becomes the following Equation (54.1). Now, the activated constraints in $\mathcal{P}_u(u_1(t)) \cap \mathcal{P}_\chi(u_1(t)|\chi_3)$ are presented in Equations (54.1) - (54.5).

$$\hat{v}_\omega(t+1) = A\hat{v}_\omega(t) + Bv_1(t+1) \tag{54.1}$$

$$\hat{x}_\omega(t+1) = \hat{x}_\omega(t) + \delta\hat{v}_\omega(t+1) \tag{54.2}$$

$$v_1(t+1) = v_1(t) + \delta u_1(t) \tag{54.3}$$

$$x_1(t+1) = x_1(t) + \delta v_1(t) + \frac{\delta^2}{2}u_1(t) \tag{54.4}$$

$$s_{\omega,1}(t) = x_1(t) - \hat{x}_\omega(t) \tag{54.5}$$

Clearly, to prove $\mathcal{P}_u(u_1(t)) \cap \mathcal{P}_\chi(u_1(t)|\chi_3) \neq \emptyset$ , we want to show there exists $u_1(t-1) \in \mathbb{U}\big(u_1(t-1)\big)$ such that the spacing at step $t$ satisfies $s_{\omega,1}(t) \geq 0$. Our poof below first demonstrates that if Condition (iv): $v_1(t-1) - \hat{v}_\omega(t-1) \geq \frac{\tilde{a}_1(t-1)\delta(B-\frac{1}{2})}{A}$ in **Theorem 1** holds at step $t-1$, then it does at step $t$ also, i.e., $v_1(t) - \hat{v}_\omega(t) \geq \frac{\tilde{a}_1(t-1)\delta(B-\frac{1}{2})}{A}$.

To begin with, we write the relative speed, $v_1(t) - \hat{v}_\omega(t)$ between $n$-CAV ($\omega$) and $s$-CAV(1) in (55.1) derived from (54.1) and (54.3).

$$\begin{aligned} v_1(t) - \hat{v}_\omega(t) &= v_1(t) - A\hat{v}_\omega(t-1) - Bv_1(t) \\ &= A[v_1(t-1) - \hat{v}_\omega(t-1)] + A\delta u_1(t-1) \end{aligned} \tag{55.1}$$

Clearly, $v_1(t-1) - \hat{v}_\omega(t-1)$ in (55.1) is bounded according to Condition (iv). We next discuss the lower bounds of the $u_1(t-1)$ in (55.1). To do that, we define a deceleration lower bound $\underline{\hat{a}}_1(t-1) = \frac{B(B-\frac{1}{2})\tilde{a}_1(t-1)}{A^2}$. Given $A = 0.1836$ and $B = 0.8163$, we have $\underline{\hat{a}}_1(t-1) = 7.66\,\tilde{a}_1(t-1) \leq \tilde{a}_1(t-1) \leq 0$. Then, given $\tilde{a}_1(t-1) = max\left[\underline{a_x}, \underline{a_{i,v}}(t-1)\right]$ in Condition (iv), we have $\underline{\hat{a}}_1(t-1) \leq max\left[\underline{a_x}, \underline{a_{i,v}}(t-1)\right]$. Recall that $\mathbb{U} = [max\left\{\underline{a_x}, \underline{a_{1,v}}(t)\right\}, min\{\overline{a_x}, \overline{a_{1,v}}(t), \overline{a_{1,d}}(t)\}]$ given in (51), we can find a $u_1(t-1) \in \mathbb{U}$, so that $u_1(t-1) \geq \underline{\hat{a}}_1(t-1)$. Then, we can further process (55.1) as follows.

$$\begin{aligned} v_1(t) - \hat{v}_\omega(t) &= A[v_1(t-1) - \hat{v}_\omega(t-1)] + A\delta u_1(t-1) \\ &\geq -\tilde{a}_1(t-1)\left(\frac{1}{2} - B\right)\delta + A\delta u_1(t-1) \\ &\geq \tilde{a}_1(t-1)\delta(B - \frac{1}{2}) + \frac{B(B-\frac{1}{2})\tilde{a}_1(t-1)\delta}{A} \\ &= \frac{\tilde{a}_1(t-1)\delta(B-\frac{1}{2})}{A} \end{aligned} \tag{55.2}$$

Built upon Condition (iv), we want to prove that there exists a control input $u_1(t-1) \in \mathbb{U}\big(u_1(t-1)\big)$ such that $s_{\omega,1}(t) \geq 0$. To show this result, we rewrite the spacing $s_{\omega,1}(t)$ in a recursive form in (56) based on (54.1) - (54.5).

$$s_{\omega,1}(t) = s_{\omega,1}(t-1) + A\delta[v_i(t-1) - \hat{v}_\omega(t-1)] + (\frac{1}{2} - B)\delta^2 u_1(t-1) \tag{56}$$

Because $\mathcal{P}(u_1(t-1), \chi)$ is feasible, $s_{\omega,1}(t-1) \geq 0$. Then $s_{\omega,1}(t) \geq A\delta[v_1(t-1) - \hat{v}_\omega(t-1)] + \left(\frac{1}{2} - B\right)\delta^2 u_1(t-1)$. According to (51), $\mathbb{U}\big(u_1(t-1)\big) = \left[max\left\{\underline{a_x}, \underline{a_{1,v}}(t-1)\right\}, min\{\overline{a_x}, \overline{a_{1,v}}(t-1), \overline{a_{1,d}}(t-1)\}\right]$. Therefore, we have $u_1(t-1) \geq max\left\{\underline{a_x}, \underline{a_{1,v}}(t-1)\right\} = \tilde{a}_1(t-1)$. Combining with the results in (55.2), we have $s_{\omega,1}(t) \geq A\delta[v_1(t-1) - \hat{v}_\omega(t-1)] + \left(\frac{1}{2} - B\right)\delta^2 u_1(t-1) \geq \tilde{a}_1(t-1)\delta^2\left(B - \frac{1}{2}\right) - \left(B - \frac{1}{2}\right)\delta^2\tilde{a}_1(t-1) = 0$. Therefore, there exists $u_1(t-1) \in \mathbb{U}\big(u_1(t-1)\big)$ so that $s_{\omega,1}(t) \geq 0$.

Furthermore, we need to consider a scenario/state transition under Scenario II. Specifically, if $s$-CAV(2) cuts in between $s$-CAV(1) and its follower, $n$-CAV($\omega$) at step $t$, Then, Scenario II transfers to the special one where $s$-CAV(2) follows $s$-CAV(1) and is followed by the $n$-CAV($\omega$). Accordingly, we need to prove the recursive feasibility of the constraints regarding $s$-CAV(2) in this special scenario. Mainly, we want to show $\mathcal{P}_u(u_2(t)) \cap \mathcal{P}_\chi(u_2(t)|\chi_1) = \mathbb{U}\big(u_2(t)\big) \neq \emptyset$, and then there exists $u_2(t) \in \mathbb{U}\big(u_2(t)\big)$ such that $\mathcal{P}_u(u_2(t)) \cap \mathcal{P}_\chi(u_2(t)|\chi_3) \neq \emptyset$ given Conditions (i) – (iv) hold. By the same arguments about Case (1), we can show that $\mathbb{U}\big(u_2(t)\big) \neq \emptyset$. Next, as the cutting-in is carried out at step $t$, we have $\eta^{l^*}_{\omega,2}(t) = 1$ and $\eta^{l^*}_{\omega,2}(t-1) = 0$ (i.e., $n$-CAV($\omega$) starts to follow $s$-CAV(2) at step $t$ but not $t-1$). Accordingly, the safety constraint (27) makes the cut-in speed $v_2(t)$ satisfy $v_2(t) - \hat{v}_\omega(t) \geq \frac{\tilde{a}_2\delta(B-\frac{1}{2})}{A}$, which indicates Condition (iv) holds. Then following the similar discussion for $s$-CAV(1) under Scenario II, we can show that $s_{\omega,2}(t+1) \geq 0$ for any control input $u_2(t) \in \mathbb{U}$. The above arguments conclude that if $s$-CAV(1) is under Scenario II, the feasibility holds after the state transition. Wrapping up the derivations and arguments above, we complete the proof for **Theorem 1**. ■

**Remark 1: Theorem 1** specifies the applicability of our MPC system. It indicates that the recursive feasibility only holds if the reaction time of HDVs satisfies $\tau_\omega \leq \bar{\tau}_\omega = \frac{-(\underline{v_x}-5\underline{a_x})-\sqrt{(\underline{v_x}-5\underline{a_x})^2-24\underline{a_x}d_\omega}}{2\underline{a_x}}$. Based on this mathematical formulation, we notice that $\bar{\tau}_\omega$ is closely related to the maximum deceleration of $s$-CAVs, $\underline{a_x}$, the speed limit $\underline{v_x}$ and Newell's displacement $d_\omega$. According to (Manual, 2000; Meng et al., 2021), we set the parameter values as follows, $\underline{v_x} = 5\ m/s$, $d_\omega = 5$, and $\underline{a_x} = -6\ m/s^2$, and then get $\bar{\tau}_\omega = 6.59\ s$. Existing literature (Manual, 2000; Wei et al., 2017) shows that the normal operator and braking system reaction time, $\tau_\omega = 2.0\ s\ \leq \bar{\tau}_\omega = 6.59$. Then, we claim that our MPC system will keep the recursive feasibility in most traffic scenarios. Furthermore, based on the CACC scheme of $n$-CAVs in (25), $A$ and $B$ can take different values depending on the value of control gain $k_2$. Our discussion of Condition (iv) above (55.2) indicates that the recursive feasibility only holds if $B > 0.5$ (or $k_2 > 0.7\ s^{-1}$ when control time gap $t_d = 0.6\ s$). This analysis result suggests the range for the CACC control parameter within which our MPC is applicable. According to Xiao et al., (2017), $k_2$ is set as $1.6\ s^{-1}$ to well fit the field data under the normal traffic and indicates our MPC system will hold the recursive feasibility in the normal traffic.

**Remark 2:** The sequential feasibility proof is closely related to the car-following models that are used to capture the moving dynamics of neighbor vehicles, i.e., HDVs and $n$-CAVs. Note that we do not control the movements of neighbor vehicles but only adopt the first-order car-following models to predict neighbor vehicles' movements. Other first-order car-following models also fit into our MPC modeling and feasibility proof. Our MPC modeling can accommodate higher-order car-following models to improve prediction accuracy but they will also increase the complexity of the feasibility proof. However, we think it is not necessary to include the higher-order car-following models for HDVs since the MPC repeatedly generates control law that can correct prediction errors resulting from previous steps. Moreover, the model complexity and significant computation expenses also degrade the benefits of adopting the higher-order models. For $n$-CAVs, we just need to use the control scheme they actually adopt in our MPC modeling. If it is a higher-order control scheme, it will complicate the recursive formulation of spacing in (56) and may give a different safety lane-change condition from Condition (iv). But we can still follow the same structure of the proof to examine the recursive feasibility.

## 7. Experiments

This section conducts three sets of numerical experiments to validate the efficiency of the synchronization control, demonstrate the merits of the adaptive weighting strategy, and investigate the effects of traffic congestion level and CAV penetration rate on the performance of the synchronization control.

The experiments are built upon the NGSIM data collected on the eastbound of the I-80 in San Francisco Bay at Emeryville, California, 4:00 pm to 4:15 pm on April 15, 2005. The testbed provides vehicle trajectory information, including lane membership, position, speed, and acceleration. We randomly select 3-lanes road segments in the testbed to generate testing cases. Among the vehicles in each testing case, we randomly choose two vehicles as the subject CAVs and label the other vehicles as neighbor vehicles. Meanwhile, we test the CAV penetration varies from 0% to 100% with 10% as step size. Correspondingly, a proportion of the neighbor vehicles are considered as neighbor CAVs for each experiment based on the predetermined CAV penetration rate. All the experiments share the default setting for some parameters presented in Table 2 according to NGSIM field data and (Police Radar Information Center, 2020). All experiments are run on a DELL Precision 3630 Tower with 3.60GHz of Intel Core i9-9900k CPU, 8 cores and 16 GB RAM in a Windows environment.

**Table 2**

Parameter settings for experiments

| Parameters | Values |
|---|---|
| MPC prediction horizon $T$ | 5 (sec) |
| Sample time interval $\delta$ | 1 (sec) |
| Vehicle length $L_i$ | 5 (m) |
| Maximum speed $\overline{v_x}$ | 33.33 (m/s) |
| Minimum speed $\underline{v_x}$ | 10 (m/s) |
| Acceleration limit $\overline{a_x}$ | 8 (m/s$^2$) |
| Deceleration limit $\underline{a_x}$ | -6 (m/s$^2$) |
| Jam density $k_c^{l,max}$ | 0.12 (veh/m) |
| Cell capacity $q_c^{l,max}$ | 2000 (veh/h) |
| Maximum weighting bound, $\xi_{max}^{I}$ ($\xi_{max}^{II}$) | 10 |

### 7.1 Performance of synchronization control with adaptive weighting strategies

We first demonstrate the performance of the synchronization control with a particular interest in showing the merit of the adaptive weighting strategy. Before presenting our experiment results, we first explain the measurement metrics used in this experiment. It is slightly different from the measurements we used in the mathematical model but provides a consistent evaluation of the performance. Recall the MPC in (36) measures the traffic performance by (i) traffic smoothness; (ii) $s$-CAVs' speeds; (iii) upstream traffic propagation. Accordingly, this experiment collects the *average traffic speed* involving both $s$-CAVs' neighbor vehicles and the upstream traffic under the synchronization control. It consistently reflects the traffic efficiency in terms of items (ii) and (iii) in the MPC. The synchronization control's impact on traffic smoothness is examined by the fluctuation of control inputs of $s$-CAVs. The smoother s-CAV movements cause less traffic fluctuation. Meanwhile, the synchronization performance under the MPC in (36) is also measured from three aspects, including (i) speed harmonization; (ii) the time that it takes for $s$-CAVs to be positioned in car-following mode; (iii) the time that it takes to stabilize the $s$-CAVs' spacing. Accordingly, the experiments collect the total time to complete the synchronization (sum of (ii) and (iii)) to demonstrate the synchronization efficiency. The speed harmonization reflects how much we stabilize the speed of the two s-CAV. It can also be evaluated by the variation of the relative speed between the two $s$-CAVs.

On the other hand, to implement the experiments, we need to determine the weights of the multi-objective function of the optimizer used in the MPC. Recall that our weighting strategy adaptively tunes $q_\eta$ in MPC-I or $q_z$ in MPC-II (see Equations (41) and (43)), but keep the weights $q_u$ for traffic smoothness and $q_v$ for the speed harmonization constant since they are consistently important during the entire synchronization control. Accordingly, before showing the merits of our adaptive weighting strategy, our experiments first want to learn the proper values for $q_u$ and $q_v$. To do that, the experiemts first screen the values of $q_u$ (and $q_{\mathcal{v}}$ ) with the step size 0.1, while using adaptive weighting strategy for other weights according to the settings in Table 3.

**Figure 3** plots the control inputs and relative speed of $s$-CAVs in some randomly picked cases. Figure 3 (a) indicates that 0.1 is the proper value for $q_u$. More exactly, this setting smoothens the movements of $s$-CAVs and leads to less traffic fluctuation as compared with the setting: $q_u = 0.0$. On the other hand, further increasing $q_u$ to 0.2 and 0.3 does not significantly improve the smoothness of $s$-CAVs' movements but results in a longer time to complete synchronization. Figure 3 (b) demonstrates that $q_{\mathcal{v}} = 0.1$ works resonably well to synchronize the relative speed (less relative speed variations), when the two $s$-CAVs are already positioned in the car-following mode (under MPC-II) and seek to stabilize their spacing and relative speed. We do not see the further improvement in the speed harmonization by increasing $q_{\mathcal{v}}$ to 0.2 and 0.3.

With the hint from **Figure 3**, we put both $q_u$ and $q_{\mathcal{v}}$ equal to 0.1 to examine the effectiveness of the adaptive weighting strategy. Specifically, we compare the comprehensive performance of the synchronization control under four weighting strategies: (i) the balanced strategy (S1-bal), which considers traffic and synchronization efficiency equally important and sets $q_\eta$ and $q_{i,\mathcal{w}}$ in MPC-I (or $q_z$ and $q_{i,\mathcal{w}}$ in MPC-II) with even weights in the objective function; (ii) the synchronization dominant strategy (S2-syn), which prioritizes the synchronization over the traffic efficiency by setting a larger value for $q_\eta$ than for $q_{i,\mathcal{w}}$ (or a larger value for $q_z$ than for $q_{i,\mathcal{w}}$ in MPC-II); (iii) the traffic dominant strategy (S3-trf), which overweighs the traffic efficiency; (vi) the adaptive weighting strategy according to Equations (41) and (43), for which we further test how the tunable parameter $\alpha$ affects its performance. The values of weights for all the weighting strategies are presented in Table 3.

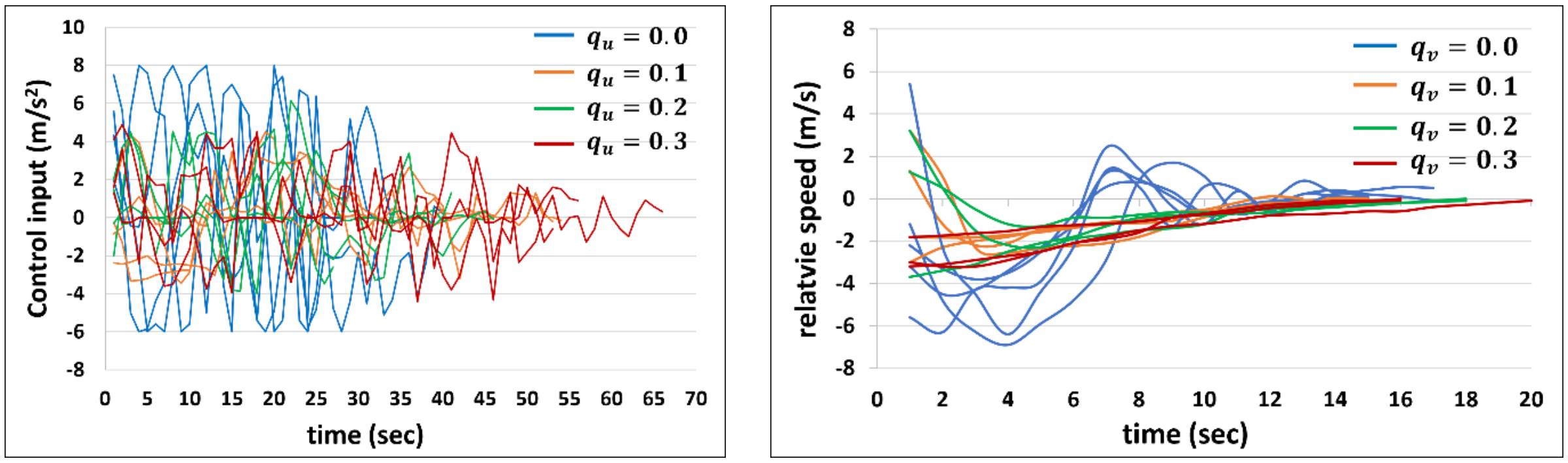

(a) Plot of control inputs of $s$-CAVs under different control weights (MPC-I & MPC-II)

(b) Plot of relative speed of $s$-CAVs under different weights (MPC-II only)

**Figure 3.** Plot of control inputs (traffic smoothness) and relative speed (speed harmonization) under different weights

**Table 3**

Values of weights for each strategy

| ***MPC-I*** | Balanced (S1) | Synchronization dominant (S2) | Traffic dominant (S3) | Adaptive weighting |
|---|---|---|---|---|
| $q_\eta$ | 0.40 | 0.40 | 0.20 | Adaptive |
| $q_{i,\mathcal{w}}$ | 0.40 | 0.20 | 0.40 | 0.40 |
| ***MPC-II*** | | | | |
| $q_z$ | 0.35 | 0.35 | 0.15 | Adaptive |
| $q_{i,\mathcal{w}}$ | 0.35 | 0.15 | 0.35 | 0.35 |

Figure 4 demonstrates the results regarding time to complete the synchronization, $s$-CAV speed and average traffic speed under the weighting strategies. It can be observed that strategy S2-syn makes the synchronization complete fastest as expected, i.e., 33.74 sec. But, it significantly slows down s-CAVs' speed (17.95 m/s on the average) and the surrounding traffic's speed (20.98 m/s on the average). On the other hand, strategy S3-trf enables $s$-CAVs and the surrounding traffic to move fastest (23.78 m/s and 21.55 m/s on average, respectively). Still, it takes the longest time to have s-CAVs complete the synchronization (113.25 sec). We can also observe that strategy S1-bal performs better than strategy S3-trf by evenly considering the time to complete the synchronization and traffic speed. More specifically, under strategy S1-bal, we can significantly reduce the synchronization time by 23.9% (from 113.25 sec under S3-trf to 86.14 sec under S1-bal) with minor sacrifice in average traffic speed (i.e., reduce 0.3% from 21.55 m/s under S3-trf to 21.49 m/s under S1-bal). However, S1-bal is not perfect. It makes the synchronization take

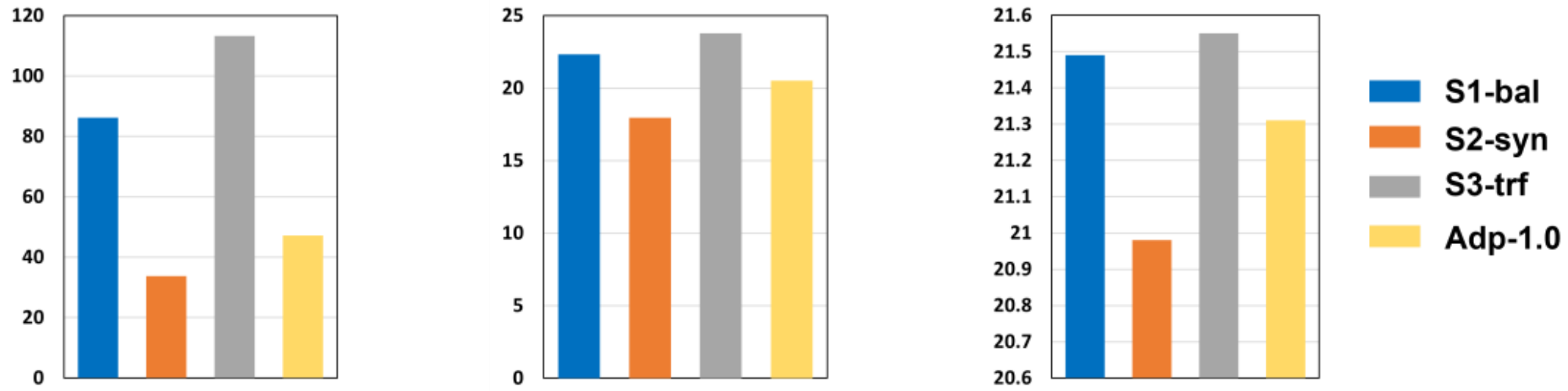

(a) Synchronization time (sec) (b) CAV speed (m/s) (c) Average traffic speed (m/s)

**Figure 4.** Comparison of traffic and synchronization efficiency under different weighting strategies.

a much longer time to complete than that under the strategy S2-syn (86.14 sec vs. 33.74 sec). On the other hand, as compared with strategy S1-bal, the adaptive weighting strategy enables the controller to trade a bit of average speed to significantly reduce the time needed to complete the synchronization. Recall that the weights $q_\eta$ in MPC-I (and $q_z$ in MPC-II) are associated with the synchronization performance metrics and they are adaptively tuned according to Equations (41) and (43), in which $\alpha$ represents the tuning level. For example, larger $\alpha$ gives larger weight on $q_\eta$ and $q_z$, which promotes the synchronization process. We tested different values of $\alpha$ and found that the performance of adaptive weighting strategy will not significantly change when $\alpha > 0.5$. As a result, we present the result under $\alpha = 1.0$ in Figure 4. The yellow column in Figure 4 indicates that under the adaptive weighting strategy (Adp-1.0), the controller trades 0.84% average speed reduction (from 21.49 m/s under S1-bal to 21.31 m/s) for a significant reduction in the time to complete the synchronization (from 86.14 sec (S1-bal) to 47.12 sec). Therefore, we conclude that the developed adaptive weighting strategy outperforms other weighting methods that have fixed weights.

### 7.2 Synchronization control effectiveness

Next, this study demonstrates the effectiveness of the synchronization control by comparing it to two field cases in NGSIM field dataset. In each case, we observed that a vehicle conducted several lane-changing and then stabilized a following mode behind the other vehicle. Even though we have no data to specify the real purpose of this movement, we consider this is a synchronization process without implementing a cooperative control, and use it as a benchmark to demonstrate the traits of our approach. Figure 5 plots the trajectories of subject vehicles (ID = 210, 202 in Case (i) and ID = 518, 516 in Case (ii)) with 1.0 second resolution from NGSIM field data (left) and under the synchronization control (right). The neighboring vehicles are labeled with corresponding IDs. To better evaluate the results, we also summarize the average speeds of neighboring vehicles and subject vehicles from field data and synchronization control in Table 4.

**Table 4**
Average vehicle speeds from field data/synchronization control

| ***Case (i)*** | | |
|---|---|---|
| Vehicle ID | Speed from field data (m/s) | Speed from MPC (m/s) |
| 202 | 33.19 | 34.81 |
| 210 | 21.73 | 20.85 |
| Neighboring vehicles | 21.22 | 22.93 |
| ***Case (ii)*** | | |
| Vehicle ID | Speed from field data (m/s) | Speed from MPC (m/s) |
| 516 | 19.85 | 21.58 |
| 518 | 26.45 | 26.50 |
| Neighbor vehicles | 27.05 | 27.28 |

More exactly, Case (i) in Figure 5 (a) provides an example under a sparse traffic only involving 5 vehicles. Figure 5 (a, left) shows that subject vehicle 210 in Case (i) initially drove on lane 1, moved from lane 1 to lane 3, and eventually approached vehicle 202 in lane 3. The process involved two lane-changing, affected three neighboring vehicles, and took 20 seconds to complete. Assuming Case (i) is an example of a catch-up maneuver, vehicles 210 and 202 are the subject CAVs. Figure 5 (a, right) provides the trajectory under the synchronization control. We can observe that vehicle 210 no longer cut in between vehicle 253 and vehicle 234 on lane 2, which very likely blocks the traffic. Instead, it stayed in lane 1 until it had a loose space to move to lane 2. After that, it stayed in lane 2, followed vehicle 253 for a while, and eventually moved to lane 3 to catch vehicle 202 when it was close. Table 4 shows that both the two subject vehicles and the neighboring vehicles drive faster under the synchronization control than they did in the field. These results demonstrate the efficiency and effectiveness of the synchronization control.

Case (i) involves a similar "catch-up" motion between vehicle 518 and vehicle 505 under congested traffic. It affects 16 vehicles within 8 seconds. The dash-dot lines in Figure 5 (b) show the trajectories of vehicle 518 observed in the field data and the simulation under the synchronization control. Even though they both involve two lane-changing maneuvers, our experiment data indicated that the synchronization took much less time (3 seconds) to complete under the synchronization control than it initially did (8 seconds). Moreover, Table 4 demonstrates a similar improvement in traffic efficiency to Case (i). These improvements in traffic and synchronization confirm the value of the complicated MPC factoring in hybrid traffic dynamics in the trajectory control. Finally, by comparing Case (ii) to Case (i), we noticed that even though Case (ii) is under congested traffic, the corresponding synchronization takes less time to complete

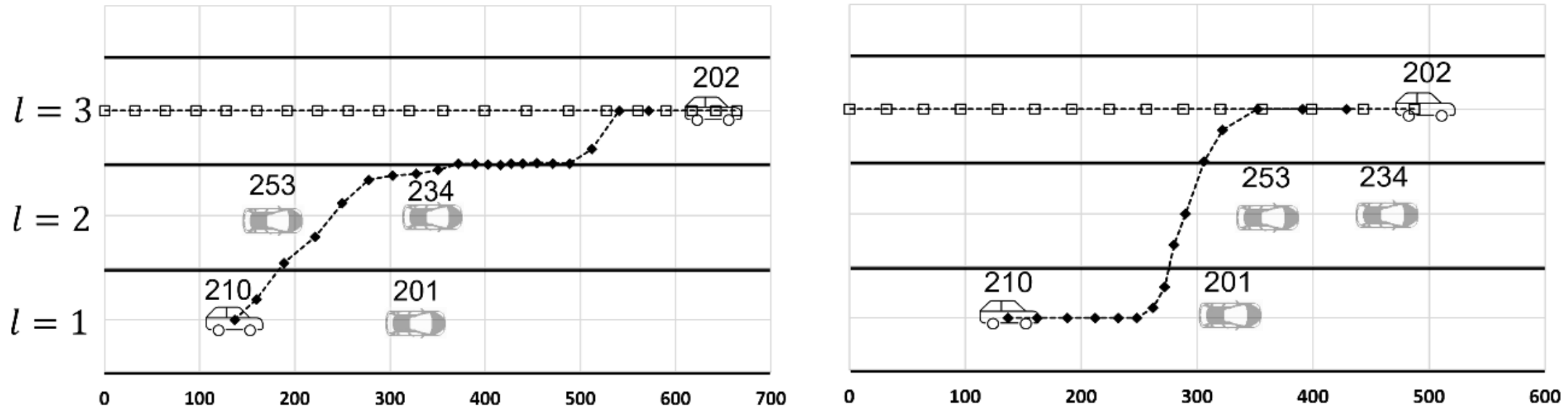

(a) vehicle trajectories in Case (i) from field data (left), and synchronization control (right).

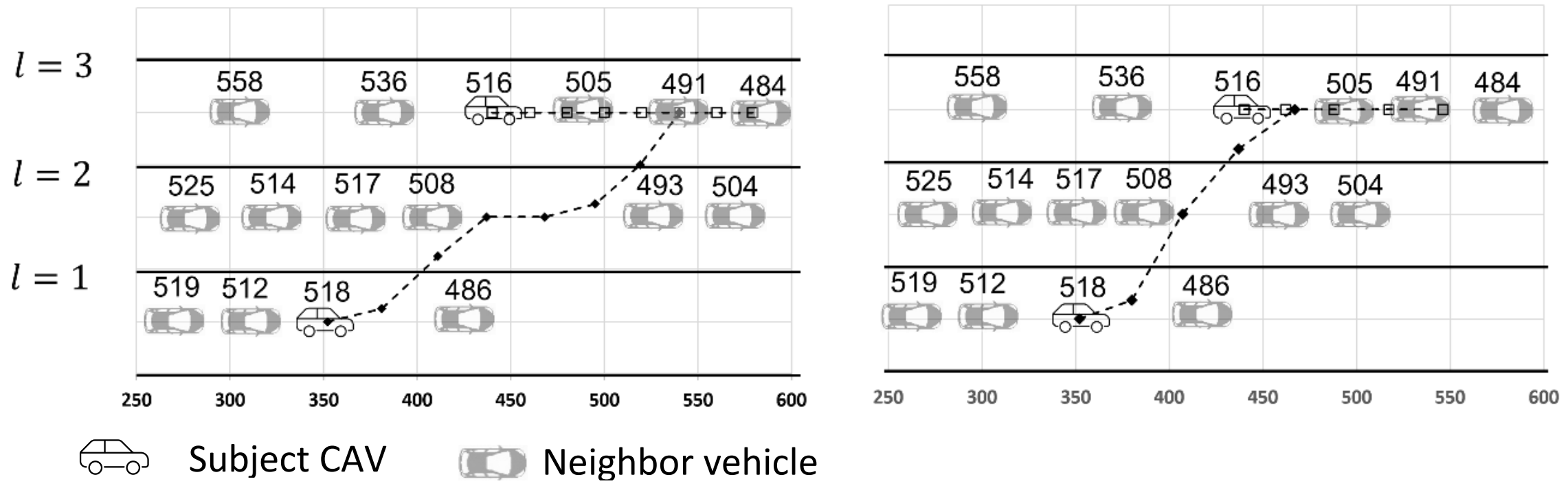

(b) vehicle trajectories in Case (ii) from field data (left), and synchronization control (right).

**Figure 5.** Comparison of vehicle trajectories from field data and by synchronization control.

than it does in Case (i) under sparse traffic. This result indicates that besides the traffic congestion level, the initial position of $s$-CAVs also plays an important role in the synchronization performance.

### 7.3 Sensitivity analysis

The results from the last section indicate that the synchronization control performance is affected by the congestion level, CAV penetration, and the initial positions of the two $s$-CAVs. This section conducts the sensitivity analysis to investigate their effects. Our experiments screened the synchronization efficiency and traffic efficiency under different traffic congestion levels represented by the traffic volume ($v$) to capacity ($c$) ratio (v/c)[2] varying from 0.35 to 0.85 and different CAV penetration ranging from 0% to 100% with the step size equal to 10%. We run the experiments in multiple cases for each scenario with a given congestion level and CAV penetration, each of which randomly selects the initial positions of the two subject CAVs. Below we discuss the experiment results.

[2] Please refer to Appendix C for a detailed discussion on the calculation of v/c.

### 7.3.1 Impact of Traffic congestion

We first test the impact of congestion levels on the synchronization control performance. Our screening experiments show that raising CAV penetration will generally benefit the synchronization control performance, but its marginal effect diminishes and becomes insensitive after 50%. To highlight the impact of traffic congestion levels, below we present the results under 50% CAV penetration in this subsection.

Figure 6 plots the time to complete the synchronization (i.e., synchronization time in the figure) and the average traffic speed of the testing cases under different congestion levels. Figure 6 (a) indicates that the time steps needed to complete the synchronization control increases gradually when the traffic v/c ratio is below 0.55 (i.e., under light and mildly congested traffic), and then rises significantly as v/c ratio increases further (i.e., under congested and highly congested traffic). This is reasonable. When the traffic is congested, the two $s$-CAVs need more time to manage their movements (e.g., lane changes, acceleration/deceleration) and then approach each other, given neighbor vehicles move slowly. On the other

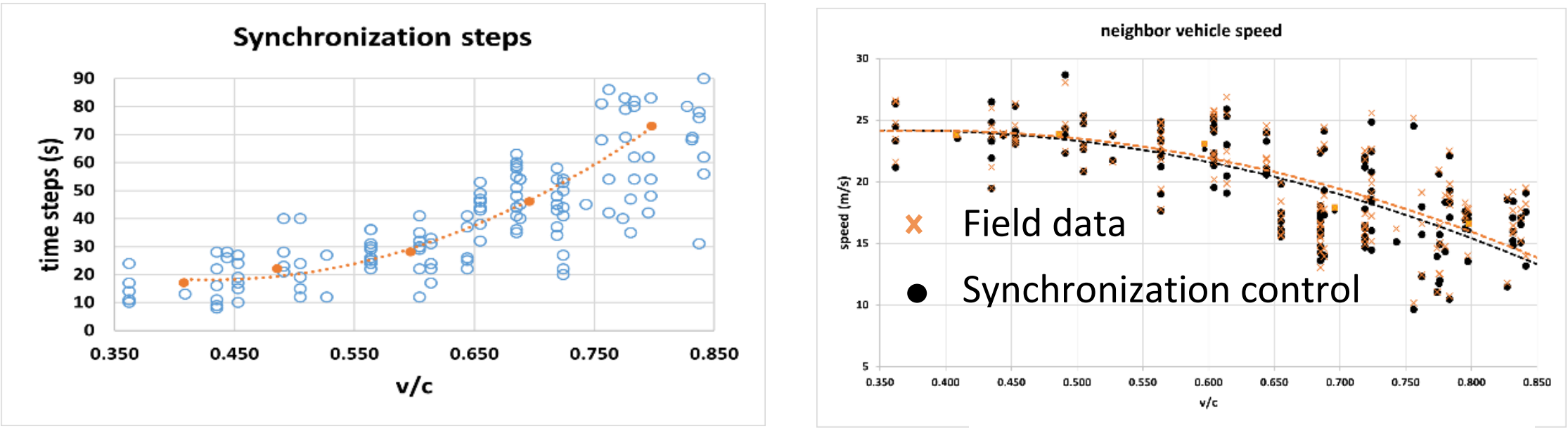

(a) synchronization time (sec)

(b) average traffic speed (m/s)

**Figure 6.** Performance of synchronization control under different congestion levels.

hand, Figure 6 (b) shows the average speed of surrounding traffic with the synchronization control in our simulation (the black curve) is very close to the corresponding field traffic speed without synchronization control (the orange curve). Therefore, we confirm that this synchronization control does not worsen traffic in most cases under different congestion levels. The negative effect only occurs slightly when traffic is highly congested such as v/c > 0.8. This observation is consistent with the findings in Section 7.2. Overall, the sensitivity analyses demonstrate that traffic congestion will slow down the synchronization process, but the synchronization control can well preserve traffic efficiency under different traffic congestion levels.

### 7.3.2 Impact of CAV penetration

Next, this study examines how CAV penetration affects the performance of the synchronization control. The above experiments have shown that raising the congestion level will generally worsen the synchronization efficiency. This effect becomes significant when v/c ratio is above 0.55. To highlight the impact of the CAV penetration, Figure 7 below presents the results of the cases with v/c varying between 0.45 to 0.55 (Note that the v/c is coming from the testing cases from field data. Under each penetration, we try to find the testing cases from the data with the v/c ratio ranging 0.45 to 0.55.)

From Figure 7 (a), we first observed that the time steps to complete the synchronization vary very much in the testing cases with different $s$-CAVs' initial positions even under the same CAV penetration, this observation renders the variation of the difficulties in synchronizing the two s-CAVs' movement. More exactly, Figure 7 (a) shows that the time needed to complete the synchronization decreases when the CAV penetration increases. We also observe that this positive effect is marginally diminishing. For instance, the mean time to complete the synchronization control is significantly reduced 27% from 38.81 sec to 28.43 sec when the CAV penetration increases from 20% to 40%. But it only reduces 11% from 28.43 sec to 25.34 sec as the CAV penetration increases from 40% to 60%. Figure 7 (b) also shows that the average traffic speed increases 1.5% from 21.75 m/s to 22.08 m/s as CAV penetration increases from 20% to 60%. Moreover, our experiments observed that s-CAVs can complete the lane change more promptly if the neighboring vehicles are mainly CAVs rather than HDVs. Therefore, we confirm that CAVs' quicker response to s-CAV movement will benefit the synchronization control performance than HDVs do. Accordingly, more CAVs present in the background traffic make the synchronization control performance better.

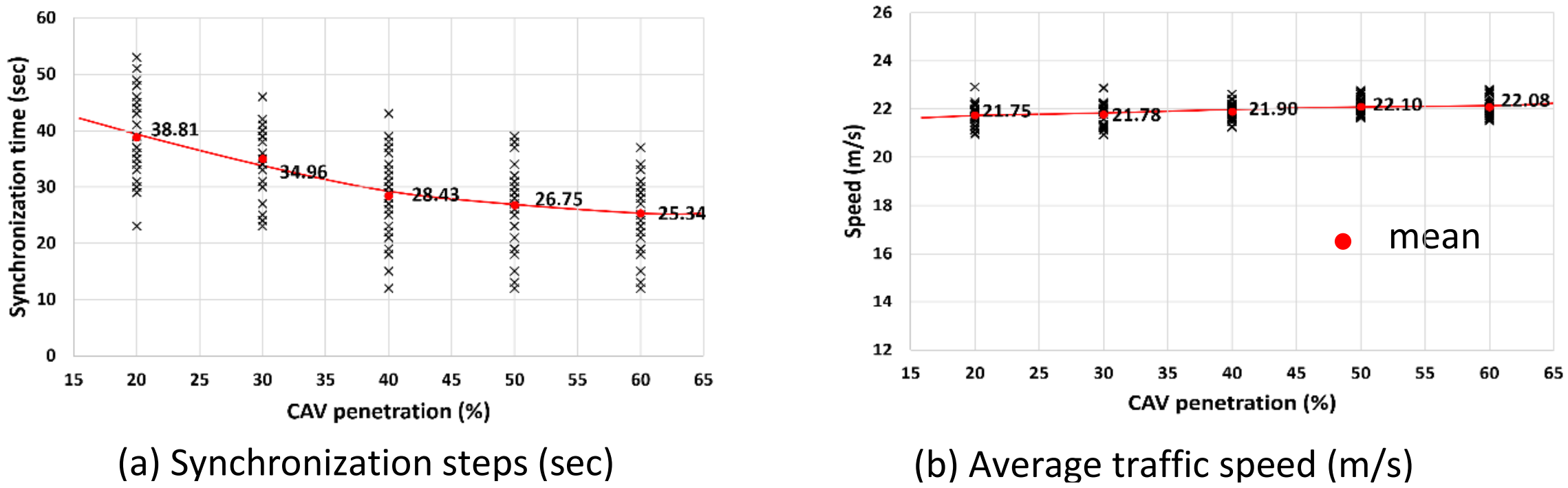

**Figure 7.** Performance of synchronization control under different CAV penetration rate.

## 8. Conclusion

This paper developed a CAV synchronization control, which cooperatively instructs the movements of two CAVs, initially apart and separated by other vehicles in a CAV and HDV mixed traffic flow, so that they can approach each other quickly and then keep a stable car-following mode. Mathematically, the synchronization control is modeled as an MPC system embedded with a mixed-integer nonlinear program (MINLP-MPC), which predicts the surrounding traffic and generates discrete optimal trajectory control laws to have the two s-CAVs catch up and synchronize their movement without jeopardizing the traffic efficiency. To couple the microscopic trajectory control and macroscopic traffic impact involved in the synchronization control, the MPC system introduces binary variables to capture the interactions between microscopic vehicle dynamics, i.e., subject CAVs and their neighbor vehicles, and macroscopic traffic propagation. We proved the recursive feasibility of the MPC model, through which we also demonstrated the applicability of the MPC system in reality. Furthermore, the optimizer of the MPC involves multiple objectives that co-consider several aspects of the traffic and synchronization efficiency. To generate a well-balanced optimal control, this study transfers the initial MPC to a hybrid MPC system, which consists of two sequential MPCs to decorrelate the objective items. Built upon that, we contribute an adaptive weighting strategy for each MPC system which enables us to adaptively tune the control priorities according to the initial traffic condition at each control time step.

The experiments based on the field data validate the effectiveness and merits of the synchronization control. Mainly, the synchronization control combined with the adaptive weighting strategy renders more space to trade a little traffic efficiency for a great improvement in the synchronization performance compared with three other strategies (S1-bal, S2-syn and S3-trf) that have fixed weights. The sensitivity analyses on the congestion levels demonstrate that the time needed to complete the synchronization stays relatively stable when traffic is light or mildly congested, i.e., v/c < 0.55, but increases significantly when traffic is highly congested, i.e., v/c > 0.75. We also examined the effect of CAV penetration. The results indicate that more CAV presenting in the traffic can help improve the synchronization efficiency, but this positive effect diminishes.

This is our first attempt to work on this synchronization control involving two CAVs. It motivates several interesting future research topics to develop this study further. For example, one of our future studies seeks to make synchronization control more resilient to communication delay and sensing error. The synchronization control depends on sensing and monitoring the neighboring vehicles' real-time movement. Even though MPC itself, as a discrete closed-loop feedback control system, can handle a certain level of prediction error and uncertainties, there are more advanced approaches to factor in those issues better. Furthermore, it is a promising research direction to develop a platoon formation controller built upon this synchronization control. Namely, we want to synchronize small platoon pairs to form a larger platoon. Those research will complicate the MPC and solution approaches developed in this study significantly, and we propose to address these new challenges in our future studies.

**Acknowledgment**: This study is partially supported by National Science Foundation, award CMMI 1818526 (Previous: 1554559)

## 10. Appendix

### Appendix A. Logical Conditions

In this appendix we reformulate logical conditions (5), (6) and (32), (33) by sets of integer constraints. Logical condition (5) models the relationship among binary variable $\rho_{j,j'}^{l}(t)$, relative positions and lane memberships of vehicle $j$ and vehicle $j'$. It can be reformulated into constraint set (57) by introducing sufficient big value $M$.

$$x_{j'}(t) - x_j(t) \le M\rho_{j,j'}^{l}(t) + M\left(1-\gamma_j^l(t)\right) + M\left(1-\gamma_{j'}^l(t)\right), \quad \forall j,j' \in J, l \in L \tag{57.1}$$

$$x_{j'}(t) - x_j(t) \ge M(\rho_{j,j'}^{l}(t) - 1), \quad \forall j,j' \in J, l \in L \tag{57.2}$$

$$\gamma_j^l(t) \ge \rho_{j,j'}^{l}(t), \gamma_{j'}^l(t) \ge \rho_{j,j'}^{l}(t), \quad \forall j,j' \in J, l \in L \tag{57.3}$$

To reformulate the logical condition (6), we introduce another binary variable $\xi_k^l(t) = \rho_{jk}^l(t)\rho_{kj'}^l(t)$, which takes value one when vehicle $j'$ is in the downstream of vehicle $j$ and there exists a vehicle $k$ between them. It can be formulated by the following constraint set (58).

$$\xi_k^l(t) \le \rho_{jk}^l(t) \quad \forall j,j' \in J, l \in L \tag{58.1}$$

$$\xi_k^l(t) \le \rho_{kj'}^l(t) \quad \forall j,j' \in J, l \in L \tag{58.2}$$

$$\xi_k^l(t) \ge \rho_{jk}^l(t) + \rho_{kj'}^l(t) - 1 \quad \forall j,j' \in J, l \in L \tag{58.3}$$

Next, we reformulate the logical condition (6) into constraint set (59), where $\delta_{j,j'}^l$ is auxiliary binary variable.

$$\delta_{j,j'}^l(t) \ge \eta_{j,j'}^l(t) \quad \forall j,j' \in J, l \in L \tag{59.1}$$

$$\rho_{j,j'}^l(t) \ge \eta_{j,j'}^l(t) \quad \forall j,j' \in J, l \in L \tag{59.2}$$

$$\eta_{j,j'}^l(t) \ge \rho_{j,j'}^l(t) + \delta_{j,j'}^l(t) - 1 \quad \forall j,j' \in J, l \in L \tag{59.3}$$

$$\sum_{k\in J/\{j,j'\}} \xi_k^l(t) \le \left(1-\delta_{j,j'}^l(t)\right)M \quad \forall j,j' \in J, l \in L \tag{59.4}$$

$$1 - \sum_{k\in J/\{j,j'\}} \xi_k^l(t) \le \delta_{j,j'}^l(t) \quad \forall j,j' \in J, l \in L \tag{59.5}$$

Finally, logical conditions (32) can be reformulated by the constraint set (60) with a sufficient large number $M$. Condition (33) can also be reformulated as similar to constraint (60.1) as in constraint (60.3).

$$\underline{x}_c^l + M\left(\phi_{i,c}^l(t) - 1\right) \le x_i(t) \le \overline{x}_c^l + M\left(1 - \phi_{i,c}^l(t)\right), \quad \forall i \in I, l \in L, c \in C_l \tag{60.1}$$

$$\underline{y}_c^l + M\left(\phi_{i,c}^l(t) - 1\right) \le y_i(t) \le \overline{y}_c^l + M\left(1 - \phi_{i,c}^l(t)\right), \quad \forall i \in I, l \in L, c \in C_l \tag{60.2}$$

$$\underline{x}_c^l + M\left(\phi_{\omega,c}^l(t) - 1\right) \le \hat{x}_\omega(t) \le \overline{x}_c^l + M\left(1 - \phi_{\omega,c}^l(t)\right), \quad \forall i \in I, l \in L, c \in C_l \tag{60.3}$$

### Appendix B. Minimum Spacing of Scenario I

In this Appendix, we derive the spacing $s_{\omega,1}(t+nk)$ according to the Equation set (52). The results are presented in Table 5, where $\Delta s_{\omega,i}(t)$ is the spacing difference of current spacing to the previous spacing, i.e., $\Delta s_{\omega,i}(t+nk) = s_{\omega,i}(t+nk) - s_{\omega,i}(t+(n-1)k)$. The results indicate that the minimum spacing is at $n = 5$.

**Table 5**
Spacing between HDV($\omega$) and $s$-CAV(1)

| ***Time*** $(\boldsymbol{t})$ | **0** | $\boldsymbol{k}$ | $\boldsymbol{2k}$ | $\boldsymbol{3k}$ |
|---|---|---|---|---|
| $x_i(t)$ | $s_{\omega,i}$ | $s_{\omega,i} + v_i\tau_\omega + \frac{1}{2}u_1(t){\tau_\omega}^2$ | $s_{\omega,i} + 2v_i\tau_\omega + 2u_1(t){\tau_\omega}^2$ | $s_{\omega,i} + 3v_i\tau_\omega + \frac{9}{2}u_1(t){\tau_\omega}^2$ |
| $\hat{x}_\omega(t)$ | 0 | $v_i\tau_\omega$ | $v_i\tau_\omega + s_{\omega,i} - d_\omega$ | $v_i\tau_\omega + 2s_{\omega,i} - 2d_\omega + \frac{1}{2}u_1(t){\tau_\omega}^2$ |
| $v_i(t)$ | $v_i$ | $v_i + u_1(t)\tau_\omega$ | $v_i + 2u_1(t)\tau_\omega$ | $v_i + 3u_1(t)\tau_\omega$ |
| $\hat{v}_\omega(t)$ | $v_i$ | $(s_{\omega,i} - d_\omega)/\tau_\omega$ | $(s_{\omega,i} - d_\omega + \frac{1}{2}u_1(t){\tau_\omega}^2)/\tau_\omega$ | $(v_i\tau_\omega + 2u_1(t){\tau_\omega}^2)/\tau_\omega$ |
| $s_{\omega,i}(t)$ | $s_{\omega,i}$ | $s_{\omega,i} + \frac{1}{2}u_1(t){\tau_\omega}^2$ | $v_i\tau_\omega + 2u_1(t){\tau_\omega}^2 + d_\omega$ | $2v_i\tau_\omega - s_{\omega,i} + 4u_1(t){\tau_\omega}^2 + 2d_\omega$ |
| $\Delta s_{\omega,i}(t)$ | | $-$ | $-$ | $-$ |

| ***Time*** $(\boldsymbol{t})$ | $\boldsymbol{4k}$ | $\boldsymbol{5k}$ | $\boldsymbol{6k}$ |
|---|---|---|---|
| $x_i(t)$ | $s_{\omega,i} + 4v_i\tau_\omega + 8u_1(t){\tau_\omega}^2$ | $s_{\omega,i} + 5v_i\tau_\omega + \frac{25}{2}u_1(t){\tau_\omega}^2$ | $s_{\omega,i} + 6v_i\tau_\omega + 18u_1(t){\tau_\omega}^2$ |
| $\hat{x}_\omega(t)$ | $2v_i\tau_\omega + 2s_{\omega,i} - 2d_\omega + \frac{5}{2}u_1(t){\tau_\omega}^2$ | $4v_i\tau_\omega + s_{\omega,i} - d_\omega + \frac{13}{2}u_1(t){\tau_\omega}^2$ | $6v_i\tau_\omega + 12u_1(t){\tau_\omega}^2$ |
| $v_i(t)$ | $v_i + 4u_1(t)\tau_\omega$ | $v_i + 5u_1(t)\tau_\omega$ | $v_i + 6u_1(t)\tau_\omega$ |
| $\hat{v}_\omega(t)$ | $(2v_i\tau_\omega - s_{\omega,i} + d_\omega + 4u_1(t){\tau_\omega}^2/\tau_\omega$ | $(2v_i\tau_\omega - s_{\omega,i} + d_\omega + \frac{11}{2}u_1(t){\tau_\omega}^2)/\tau_\omega$ | $(v_i\tau_\omega + 6u_1(t){\tau_\omega}^2)/\tau_\omega$ |
| $s_{\omega,i}(t).$ | $2v_i\tau_\omega - s_{\omega,i} + \frac{11}{2}u_1(t){\tau_\omega}^2 + 2d_\omega$ | $v_i\tau_\omega + 6u_1(t){\tau_\omega}^2 + d_\omega$ | $s_{\omega,i} + 6u_1(t){\tau_\omega}^2$ |
| $\Delta s_{\omega,i}(t)$ | $-$ | $-$ | $+$ |

## Appendix C. Calculation of v/c by Greenshield Model

In this Appendix, we introduce how to calculate the v/c based on the traffic information of testing cases. Recall that we randomly select road segments in the NGSIM dataset to generate testing cases. Each testing case contains the information of flow density, vehicle spacings and vehicle velocities. Based on this information, we estimate the congestion level (i.e., the traffic volume ($v$) to capacity ($c$) ratio) of each testing case through Greenshield model in Equation (61), where $k_j$ is the jam density and $k$ is the density.

$$\frac{v}{c} = \frac{4c}{k_j}k - \frac{4c}{k_j^2}k^2 \tag{61}$$

The jam density $k_j$ and capacity c adopt the values in Table 2.